\documentclass[12pt]{amsart}
\usepackage[letterpaper,margin=1in,footskip=0.25in]{geometry}

\usepackage{microtype}
\usepackage{eucal}
\usepackage{setspace}
\usepackage{mathrsfs}
\usepackage{quiver}

\usepackage[utf8x]{inputenc}

\usepackage{amsmath,amssymb,amsthm}
\usepackage{mathtools}
\makeatletter
  \newcommand{\supsize}{%
    \expandafter\ifx\csname S@\f@size\endcsname\relax
      \calculate@math@sizes
    \fi
    \csname S@\f@size\endcsname
    \fontsize\sf@size\z@\selectfont
  }
  \DeclareRobustCommand{\tsup}[1]{%
    \leavevmode\raise.9ex\hbox{\supsize #1}%
  }
  \DeclareTextSymbolDefault{\textprimechar}{OMS}
  \DeclareTextSymbol{\textprimechar}{OMS}{48}
  \DeclareRobustCommand{\tprime}{\tsup{\textprimechar}}
  \ProvideTextCommandDefault{\cprime}{\tprime}
\makeatother

\usepackage{mathabx}

\usepackage[shortlabels]{enumitem}
\setlist{noitemsep}

\usepackage[pdfusetitle,colorlinks]{hyperref}
\hypersetup{allcolors=blue}
\usepackage[capitalise,noabbrev]{cleveref}

\usepackage{colonequals}

\usepackage{graphicx}
\usepackage{extarrows}
\usepackage{tikz,tikz-cd}

\usepackage[french, english]{babel}

\crefformat{equation}{\ensuremath{(#2#1#3)}}
\crefmultiformat{equation}{\ensuremath{(#2#1#3)}}{ and~\ensuremath{(#2#1#3)}}{, \ensuremath{(#2#1#3)}}{, and~\ensuremath{(#2#1#3)}}

\numberwithin{figure}{section}
\numberwithin{equation}{section}
\theoremstyle{definition}

\newtheorem{lemma}[figure]{Lemma}

\newtheorem{question}[figure]{Question}

\newtheorem{thm}[figure]{Theorem}
\newtheorem{lem}[figure]{Lemma}

\newtheorem{cor}[figure]{Corollary}
\newtheorem{prop}[figure]{Proposition}
\newtheorem{defn}[figure]{Definition}

\newtheorem{rmk}[figure]{Remark}
\newtheorem{exam}[figure]{Example}
\newtheorem{quest}[figure]{Question}

\newtheorem{claim}[figure]{Claim}

\newtheoremstyle{cited}{.5\baselineskip\@plus.2\baselineskip\@minus.2\baselineskip}{.5\baselineskip\@plus.2\baselineskip\@minus.2\baselineskip}{\itshape}{}{\bfseries}{\bfseries .}{5pt plus 1pt minus 1pt}{\thmname{#1}\thmnumber{ #2}\thmnote{ \normalfont #3}}
\theoremstyle{cited}

\newcommand{\an}{\operatorname{an}}
\newcommand{\Ann}{\operatorname{Ann}}
\newcommand{\height}{\operatorname{ht}}

\newcommand{\perfd}{\operatorname{perfd}}
\newcommand{\rk}{\operatorname{rk}}

\newcommand{\spc}{\operatorname{sp}}

\newcommand{\trdeg}{\operatorname{tr.deg}}

\newcommand{\Spec}{\operatorname{Spec}}
\newcommand{\Spf}{\operatorname{Spf}}

\newcommand{\Min}{\operatorname{Min}}
\newcommand{\Spa}{\operatorname{Spa}}

\newcommand{\ins}{\subset}

\newcommand{\Image}{\operatorname{Im}}

\newcommand{\surj}{\twoheadrightarrow}

\newcommand{\Zar}{\operatorname{Zar}}

\newcommand{\coperf}{\operatorname{coperf}}

\newcommand{\argnorm}{\operatorname{argnorm}}
\newcommand{\res}{\operatorname{res}}

\begin{document}

\title{Abhyankar valuations, Prüfer-Manis valuations, and perfectoid Tate algebras}
\author{Dimitri Dine, Jack J Garzella}
\maketitle

\begin{abstract}
    Let \(K\) be a perfectoid field. We describe all quotient fields of the perfectoid Tate algebra
    \begin{equation*}T_{n, K}^{\perfd}=K\langle X_{1}^{1/p^{\infty}},\dots,
        X_{n}^{1/p^{\infty}}\rangle
    \end{equation*}in any number $n\geq1$ of variables in terms of (completed perfections of)
    the nonarchimedean fields
    $K_{r_1,\dots, r_l}$ occuring in Berkovich geometry. We prove that every quotient field
    \begin{equation*}L=T_{n, K}^{\perfd}/\mathfrak{m}\end{equation*}is a so-called \textit{semi-immediate}
    extension of $K_{r_1,\dots, r_l}^{\perfd}$ for some
    \begin{equation*}l\leq \min(n-\height(\mathfrak{m}^{\flat}\cap (T_{n, K^{\flat}})^{\coperf}), n-1),
    \end{equation*}which pins down the value groups and the residue fields of the possible quotient fields $L$. 
    Moreover, we show that if
    \begin{equation*}\mathfrak{m}^{\flat}\cap (T_{n, K^{\flat}})^{\coperf}\neq 0,\end{equation*}
    at least one of the radii $r_{i}$
    has to be irrational, i.e., 
    \begin{equation*}r_{i}\not\in \sqrt{\vert K^{\times}\vert}.\end{equation*}
    The main ingredient in our proof is the notion of \textit{topologically simple} valuations, which generalize
    type (IV) points
    in the classification of points on $\Spa(K\langle T\rangle)$ to the case of higher-dimensional polydisks. 
    We also consider \textit{rational Abhyankar} valuations and \textit{irrational Abhyankar} valuations,
    which generalize type (II) and (III) points, respectively. 
    We deduce our main result 
    from a description of topologically simple absolute values and of Abhyankar absolute values
    on the usual Tate algebra.
    Along the way, we also show that our topologically simple valuations are the same as
    Pr\"ufer-Manis valuations of Knebusch-Zhang.
    
    Finally, we also show that all allowed possibilities for the quotient fields $L$
    do indeed occur (i.e., the above bound $l\leq n-1$ is optimal) by
    generalizing an example of Gleason.
\end{abstract}

\tableofcontents

\section{Introduction}\label{sec:introduction}

Over the past 60 or so years,
Tate algebras have proven very useful in the study of \(p\)-adic geometry.
Remarkably, despite not being essentially of finite type over their
base field, they give a good theory of geometry which is has similar
properties to finite type geometry.\footnote{
    rings appearing in such situations are called
    \textit{topologically of finite type} in the literature.
}
At this point, we know proofs of (or analogues of, or alternatives to)
many of the fundamental theorems of algebraic 
geometry, as may be found in e.g. \cite{Atiyah-MacDonald}:
Nullstellensatz, noether normalization, dimension theory, and so on.

Since their introduction in 2010, perfectoid Tate rings have been important objects in \(p\)-adic
geometry. 
For example, they are important for various results in the Langlands program
\cite{Fargues-Scholze, scholze-2015-torsion, caraiani-scholze-2017-shimura-1}. The prototypical example of a perfectoid Tate ring which is not a field is the perfectoid Tate algebra in $n\geq1$ variables\begin{equation*}T_{n, K}^{\perfd}=K\langle X_{1}^{1/p^\infty},\dots, X_{n}^{1/p^{\infty}}\rangle,\end{equation*}whose quotients can be viewed as analogs in perfectoid theory of classical affinoid algebras in the sense of Tate. As perfectoid Tate algebras have become better studied, there is an 
important meta-question about them which remains unresolved.

\begin{quest}
    \label{quest:perfd:good}
    Do quotients of perfectoid Tate algebras $T_{n, K}^{\perfd}$, for $n\geq1$, give a good theory of 
    ``perfectoid finite type'' geometry, or are they
    badly behaved objects which should only be studied on functors of points?
\end{quest}

The question is not completely well-defined, 
but the heart of the question is commutative algebra. 
A positive answer to Question \ref{quest:perfd:good} would consist
of proofs of (or analogues of, or alternatives to) 
the various commutative-algebraic theorems that form
the basis of
algebraic/analytic geometry, as might be found in 
\cite{Atiyah-MacDonald} or \cite{BGR}.
A negative answer would be some sort of counterexample showing
that such theorems can never hold.
Of course, Question \ref{quest:perfd:good} is ill-defined as stated,
and is open to interpretation.
Nonetheless, whatever precise answers one can get will likely
be important for the applications of perfectoid spaces.

There is some evidence for a positive answer to 
Question \ref{quest:perfd:good},
for example the equivalence of ``Zariski closed''
and ``strongly Zariski closed'' subsets from \cite[Remark~7.5]{Prismatic}, 
or Gleason's result in \cite{Gleason22N} proving the Nullstellensatz
in the case when \(K\) is algebraically closed
and spherically complete with value group $\mathbb{R}_{>0}$ and uncountable residue field.

However, other recent work might suggest a negative answer
to Question \ref{quest:perfd:good}, for example
the apparent hardness of the dimension problem 
(see \cite[Problem~I.11.1]{Fargues-Scholze}, 
\cite{garzella-2024-perfectoid-uncountable}).
Another work, which shall be central to our own story,
is \cite{Gleason22N}, which
shows that the perfectoid Tate algebra in three\footnote{
    The paper also shows this for the perfectoid
    Tate algebra in two variables, with slightly
    less control on the quotient field.
}
or more variables surjects onto a field \(C(y)\)
which is not a finite extension of the base field \(K\).
Thus, perfectoid Tate algebras do not satisfy 
the Nullstellensatz that polynomial rings and
classical Tate algebras satisfy,
at least in the case when the (divisible closure of the) value group $\vert K^{\times}\vert$ of \(K\) is
not all of \(\mathbb{R}_{>0}\).

However, as Gleason observes, the field \(C(y)\) is precisely
the completed perfection $K_{r}^{\perfd}$ of the completed residue field $K_{r}$ of a type \(III\) point in Berkovich geometry. The fields \(K_{r}\),
and their higher analogues \(K_{r_{1}, \ldots, r_{n}}\),
appear in an analogue of the Nullstellensatz
for affinoid algebras in the sense of Berkovich
(e.g. \(K\langle T / r\rangle \) for some radius \(r \notin \sqrt{\vert K^{\times}\vert}\)), see Ducros \cite{Ducros07}, Théorème 2.7.

It turns out that Gleason's example, and the appearance of Berkovich's field $K_{r}$ in it, can be
embedded into a general valuation-theoretic context. Note that quotient fields of the perfectoid
Tate algebra in any number of variables are \textit{``perfectly topologically finitely generated''} nonarchimedean field extensions of $K$, i.e., each of them has a dense subfield generated by the $p$-power roots of a finite set of elements. It can be shown using the theory of Abhyankar valuations (see Corollary~\ref{Temkin's remark 4}, which follows from a variant of \cite{Temkin13}, Remark~2.1.2) that every such nonarchimedean extension $L/K$ of $K$ (say, of characteristic $p$) factors as \begin{equation*}K\subseteq L^{\prime}=K_{r_{1},\dots, r_{l}}^{\perfd}\subseteq L,\end{equation*}where $\vert L^{\times}\vert/\vert L^{\prime\times}\vert$ is torsion and the residue field of $L$ is algebraic over that of $L^{\prime}$ (we call such extensions $L/L'$ semi-immediate). Therefore, to describe all possible quotient fields $L$ of the perfectoid Tate algebra $T_{n, K}^{\perfd}$, it suffices to: \begin{enumerate}[(1)]\item Bound the number $l$ of radii $r_1,\dots, r_l$ that can occur in the above factorization of $L$, \item Determine the $\mathbb{Q}$-rank of the multiplicative abelian group \begin{equation*}\vert L^{\times}\vert/\vert K^{\times}\vert\end{equation*}(for example, can it occur that some/all of the radii $r_1,\dots, r_l$ belong to $\sqrt{\vert K^{\times}\vert}$?), and \item Restrict the possibilities for the semi-immediate extension $L/L'$; for example, one could ask whether this extension has to be the completed perfection of a finite extension of $K_{r_{1},\dots, r_{l}}$. \end{enumerate}
In this paper, we fully address part (1) of this problem, i.e., we find an upper bound (namely, $l\leq n-1$) on the number $l$ of radii $r_1,\dots, r_l$ for a quotient field of the perfectoid Tate algebra in any number $n\geq1$ of variables, and prove that this upper bound is optimal. We also obtain partial results on part (2) by imposing additional hypotheses either on the ground field $K$ or on the tilt of the quotient map $T_{n, K}^{\perfd}\to L$. 
\begin{thm}[Theorem \ref{Main theorem, in the body of the paper}]\label{Main theorem, introduction}   
    Let $\mathfrak{m}\subsetneq T_{n, K}^{\perfd}$ be a maximal ideal of the perfectoid Tate algebra $T_{n, K}^{\perfd}=K\langle X_{1}^{1 / p^{\infty}}, \ldots, X_{n}^{1 / p^{\infty}}\rangle$ and let \begin{equation*}L=T_{n, K}^{\perfd}/\mathfrak{m}\end{equation*}be the corresponding quotient field. Then there exists a polyradius $(r_1,\dots, r_l)\in (0, 1]^{l}$ such that $L$ is a semi-immediate extension of $K_{r_{1},\dots, r_{l}}$, where \begin{equation*}l\leq \min(n-\height(\mathfrak{m}^{\flat}\cap (T_{n, K^{\flat}})^{\coperf}), n-1),\end{equation*}where $\mathfrak{m}^{\flat}$ is the tilt of $\mathfrak{m}$ in the sense of \cite{Dine22}, \S4. Moreover, if $K$ has uncountable residue field, then none of the radii $r_1,\dots, r_l$ is rational, i.e.,\begin{equation*}r_1,\dots, r_l\not\in \sqrt{\vert K^{\times}\vert},\end{equation*}and if $K$ is arbitrary, but \begin{equation*}\mathfrak{m}^{\flat}\cap (T_{n, K^{\flat}})^{\coperf}\neq0,\end{equation*}then there exists at least one $i\in\{1,\dots, l\}$ with \begin{equation*}r_{i}\not\in \sqrt{\vert K^{\times}\vert}.\end{equation*} 
\end{thm}
The proof of Theorem \ref{Main theorem, introduction} is a bit subtle (especially, the proof of the last assertion on the irrationality of at least one $r_{i}$ when the residue field of $K$ is countable),
and relies on the theories of Abhyankar valuations and so-called Prüfer-Manis valuations. We will give an indication of the proof below.

In our second main theorem, we show that all of the possibilities
afforded by the bound $l \leq n-1$ in Theorem \ref{Main theorem, introduction} are
realized. 
\begin{thm}[Theorem \ref{thm:multivar:gleason}, Corollary \ref{cor:perftate:finext}, Corollary \ref{cor:perftate:finext2}]
    \label{thm:perftate:lowerbound}
    For every $n\geq 3$ and every $l\leq n-2$, there exists a continuous surjection \begin{equation*}T_{n, K}^{\perfd}\twoheadrightarrow K_{r_{1},\dots, r_{l}}^{\perfd}.\end{equation*}Moreover, for every $n\geq2$, there exists a continuous surjection \begin{equation*}T_{n, K}^{\perfd}\twoheadrightarrow L,\end{equation*}where $L$ is a semi-immediate extension of $K_{r_{1},\dots, r_{n-1}}^{\perfd}$ which moreover admits $K_{s_{1},\dots, s_{n-1}}^{\perfd}$, for a possibly different polyradius $(s_1,\dots, s_{n-1})$, as a finite separable extension.\end{thm}
We note that, while the above theorem shows that the bound $l\leq n-1$ in Theorem \ref{Main theorem, introduction} is optimal, we do not know whether the more precise bound \begin{equation*}l\leq \min(n-\height(\mathfrak{m}^{\flat}\cap (T_{n, K^{\flat}})^{\coperf}, n-1)\end{equation*}given in that theorem can be improved upon in the case when \begin{equation*}\mathfrak{m}^{\flat}\cap (T_{n, K^{\flat}})^{\coperf}\neq0.\end{equation*}In the case when the value group of $K$ is all of $\mathbb{R}_{>0}$ (with no assumptions on the residue field of $K$), our results on topologically simple and Abhyankar valuations yield the sharper bound \begin{equation*}
l\leq n-\height(\mathfrak{m}^{\flat}\cap (T_{n, K^{\flat}})^{\coperf})-1,\end{equation*}see
Proposition \ref{Any number of variables and large value group}. We thus leave the following
question open.
\begin{question}Let $K$ be an arbitrary perfectoid field, let $n\geq2$ and let \begin{equation*}\mathfrak{m}\subsetneq T_{n, K}^{\perfd}\end{equation*}be a maximal ideal in the perfectoid Tate algebra in $n$ variables over $K$ such that \begin{equation*}\mathfrak{m}^{\flat}\cap (T_{n, K^{\flat}})^{\coperf}\neq0.\end{equation*}In this case, can the bound on $l$ in Theorem \ref{Main theorem, introduction} be improved to \begin{equation*}l\leq n-\height(\mathfrak{m}^{\flat}\cap (T_{n, K^{\flat}})^{\coperf})-1?\end{equation*}\end{question} 

\subsection{Methods of proof}

Our proof of Theorem \ref{Main theorem, introduction} is of a valuation-theoretic nature. Its main ingredient is the notion of topologically simple valuation on an affinoid algebra $A$, which generalizes type (IV) from Berkovich's description of points on the closed unit disk to higher-dimensional polydisks, and its relationship with the notion of Abhyankar valuations. 
Let \(A\) be an affinoid algebra in the sense of Berkovich. 
A continuous rank \(1\) valuation \(v\) on \(A\) is \textit{Abhyankar} if it satisfies equality in (a continuous analogue of) Abhyankar's inequality.
Abhyankar valuations on \(K\langle T\rangle \) are precisely type (II) and (III) points. We can then differentiate Abhyankar valuations into \textit{rational} and \textit{irrational} ones, which correspond to type (II) points and type (III) points, respectively. 
Furthermore, we say that a rank $1$ valuation \(v\) on $A$ is \textit{topologically simple} if the topology on $A$ defined by the valuation has no nonzero closed ideals.
Topologically simple rank-$1$ valuations on the one-variable Tate algebra \(K\langle T\rangle \) are precisely the type (IV) points, see Proposition \ref{Topologically simple norms on the polynomial algebra}.
We show that the notion of topologically simple valuations is equivalent to the notion of Pr\"ufer-Manis valuations of Knebusch and Zhang
\cite{Knebusch-Zhang}.

Now, the statement of Theorem \ref{Main theorem, introduction} reduces by tilting to the the case when $K$ is of characteristic $p$. In this case, we show that every quotient field $L$ of the perfectoid Tate algebra $T_{n, K}^{\perfd}$ which is a counterexample to the classical Nullstellensatz for $T_{n, K}^{\perfd}$ (i.e., $L$ is not a finite extension of $K$) gives rise to a topologically simple continuous rank-$1$ valuation $v$ on a positive-dimensional affinoid algebra $A$ in the sense of Tate, namely, the restriction of the absolute value of $L$ to the image of $T_{n, K}\subseteq T_{n, K}^{\perfd}$ inside $L$. Choosing a Noether normalization \begin{equation*}T_{d, K}\hookrightarrow A,\end{equation*}we use a description of topologically simple valuations on $T_{d, K}$ (and, more precisely, the fact that a topologically simple valuation on $T_{d, K}$ is not Abhyankar) to show that the quotient field $L$ is a semi-immediate extension of $K_{r_{1},\dots, r_{l}}^{\perfd}$ for \begin{equation*}l\leq \min(d=n-\height(\mathfrak{m}\cap (T_{n, K})^{\coperf}), n-1),\end{equation*}which settles the main assertion of Theorem \ref{Main theorem, introduction}. Moreover, the statement on the irrationality of at least one of the radii $r_{i}$ (for arbitrary ground fields $K$) follows from the following valuation-theoretic result. 
\begin{thm}
    [Theorem \ref{Topologically simple rules out Abhyankar}]
    \label{thm:ab:topsimp}
    Let \(A\) be an affinoid algebra in the sense of Tate (i.e. a strictly affinoid algebra in the sense of Berkovich) which is not a field. Then there does not exist a continuous absolute value \(v\) on $A$ (over $K$) that is both topologically simple and a rational Abhyankar valuation.
\end{thm}
The proof of Theorem \ref{thm:ab:topsimp} uses the notion of the Shilov boundary for a strictly affinoid algebra $A$ and its description, due to Berkovich, in terms of minimal prime ideals of the reduction\begin{equation*}\widetilde{A}=A^{\circ}/A^{\circ\circ}.\end{equation*}Indeed, the rational Abhyankar valuations turn out to be the same as the weakly Shilov points of $\Spa(A, A^{\circ})$ in the sense of Bhatt and Hansen (see \cite{Bhatt-Hansen}, Proposition 2.9), i.e., those points which become points of the Shilov boundary on some rational neighborhood, and our proof of Theorem \ref{thm:ab:topsimp} proceeds by showing that a topologically simple valuation on a strictly affinoid algebra cannot be weakly Shilov. On the other hand, we caution the reader that there exist examples of continuous rank-$1$ valuations on affinoid algebras in the
sense of Berkovich, and even on affinoid algebras in the sense of Tate, which are simultaneously (irrational) Abhyankar and topologically simple (indeed, the field $K_{r}$, for $r\not\in \sqrt{\vert K^{\times}\vert}$,
is the affinoid algebra of the ``irrational circle'' $\{\, v\in \mathscr{M}(K\langle T\rangle)\mid v(T)=r\,\}$). 

The proof of Theorem \ref{thm:perftate:lowerbound}
closely follows the strategy of Gleason in 
\cite{Gleason22N}. 
That argument depends on a certain division algorithm for power
series fields with a certain finiteness property;
Gleason then proceeds to construct a special element whose coefficients
can be manipulated to satisfy the inputs of the division algorithm
argument. 
Finally, Gleason uses a Noether normalization statement for
completed perfections to remove
an extra variable.
We present each step in a slightly more general context than
\cite{Gleason22N}, so as to achieve the statement in 
Theorem
\ref{thm:perftate:lowerbound}.
The first and third step closely follow the proof in \cite{Gleason22N},
while the second requires a slightly more technical argument.

\subsection{Outline}

We conclude this introduction with a brief outline of the individual sections of the paper. In
Section \ref{sec:background} we collect some notation and terminology used in the paper. In Section
\ref{sec:abhyankar:val} we recall the definition of Abhyankar valuations à la \cite{Ducros07}. We also recall the equivalence between the notions of rational Abhyankar valuations and weakly Shilov points, already observed by Bhatt and Hansen in \cite{Bhatt-Hansen}, Proposition 2.9, and we include a self-contained proof of this result, since the proof technique is related to our proof of Theorem \ref{thm:ab:topsimp} later in Section \ref{sec:topologically simple}.
In Section \ref{sec:topologically simple}, we introduce topologically simple valuations,
prove a few basic facts, and characterize the topologically simple valuations
on the affine line and the affinoid ball. We give an example of a continuous valuation on an affinoid
algebra that is simultaneously topologically simple and (irrational) Abhyankar, and then we prove Theorem \ref{thm:ab:topsimp} which ensures that, at least for strictly affinoid algebras, such examples cannot occur if we require the Abhyankar valuation to be rational. 
In Section \ref{sec:proof of main theorem}, we combine our results from the previous sections and prove Theorem \ref{Main theorem, introduction}. 
Finally, in Section \ref{sec:Gleason} we modify Gleason's argument from \cite{Gleason22N} to prove
Theorem \ref{thm:perftate:lowerbound}.

\subsection{Acknowledgements}

The authors wish to thank their PhD advisor, Kiran Kedlaya, for many helpful conversations, advice,
and support. They also thank Ian Gleason for helpful comments on a previous draft of this work,
including clarifying Proposition 4.5 of \cite{Gleason22N}. 
They also thank Jared Weinstein, Kevin Ho, Shubhankar Sahai, Nathan Wenger, and Nik Castro for helpful
conversations. 
They also thank Vitor Borges and Tiklung Chan 
for inviting the second-named author to give a talk at the UCSD Junior Analysis Seminar, from which the discussions that led
to this paper started.
The first-named author would also like to thank Eva Viehmann for supervising his undergraduate thesis on the (one-variable) perfectoid Tate algebra back in 2019 and for the helpful conversations they had during the process of preparation of that thesis. The first-named author is also grateful to Kazuma Shimomoto for helpful conversations about that author's previous work on Shilov boundaries, and to Konstantin Ardakov for pointing out a mistake in the said earlier work, which helped us find and correct an error in the statement and proof of Theorem \ref{thm:ab:topsimp} in a previous version of this paper. 
The authors are also grateful to an anonymous reviewer, who pointed out an inaccuracy in a previous
version of this draft.

The second-named author was partially supported by the National Science Foundation Graduate 
Research Fellowship Program under Grant No. 2038238, and a fellowship from the Sloan Foundation.

\section{Background and Notation}\label{sec:background}

\subsection{Conventions}

In this paper, the word \textit{residue field} always denotes the residue field of a nonarchimedean
field. If we discuss residue fields in the geometric sense of Berkovich/adic spaces, we use the term
\textit{completed residue field}. Finally, if we discuss fields which are quotients of a ring, which
could perhaps be called ``algebraic residue fields'', we will always use the term \textit{quotient
field}.

In what follows, we will often need to distinguish between the perfect closure (i.e., colimit
perfection, also known as the coperfection) and the completed (co)perfection of a uniform Banach algebra in characteristic $p$.
 In the literature, both of these could be notated by \(-^{\text{perf}}\).
In order to avoid confusion, we don't use this notation.
If \(A\) is a uniform \(K\)-Banach algebra over a nonarchimedean field $K$ of characteristic $p$, then
\(A^{\coperf}\) denotes the perfect closure (coperfection)
and \(A^{\perfd}\) denotes the completed perfection, i.e., the completion of $A^{\coperf}$ with respect to the canonical extension of a power-multiplicative norm defining the topology on $A$.

We will often need to distiguish between affinoid algebras in the sense of Tate and affinoid
algebras in the sense of Berkovich. We make much use of both notions, so we try to specify wherever
possible, unless the context is clear. Thus, we always use ``affinoid in the sense of Tate'' or
``strictly affinoid'' for the former notion, and we use ``affinoid in the sense of Berkovich'' or
``Berkovich-affinoid'' for the latter.

\subsection{Normed Rings}

All rings are commutative and with $1$. We use multiplicative (not additive) notation for valuations on rings. Recall that a seminorm on an abelian group $G$ is a function $\lVert \cdot \rVert: G \to \mathbb{R}_{\geq 0}$ such that $\lVert 0 \rVert = 0$ and such that the triangle inequality $\lVert f + g \rVert \leq \lVert f \rVert + \lVert g \rVert$ holds for any two elements $f, g \in G$. A seminorm is said to be nonarchimedean if for any $f, g \in G$ the stronger inequality $\lVert f + g \rVert \leq \max(\lVert f \rVert, \lVert g \rVert)$ is satisfied. Unless explicitly stated otherwise we only consider nonarchimedean seminorms. A seminorm $\lVert \cdot \rVert$ is called a norm if $0$ is the only element of $G$ which is mapped to $0$ under $\lVert \cdot \rVert$.

For a ring $A$, a seminorm on $A$ is a seminorm on the underlying additive abelian group of $A$ such that $\lVert 1\rVert=1$. A seminorm $\lVert \cdot \rVert$ on a ring $A$ is said to be submultiplicative (or a ring seminorm) if it is compatible with multiplication in the sense that \begin{equation*} \lVert fg \rVert \leq \lVert f \rVert \lVert g \rVert \end{equation*}for any two elements $f, g \in A$. By a (semi)normed ring (respectively, a Banach ring) $(A, \lVert \cdot \rVert)$ we will mean a ring $A$ equipped with a submultiplicative (semi)norm (resp., a complete submultiplicative norm) $\lVert \cdot \rVert$. In this paper we tacitly assume all our seminorms on rings to be submultiplicative and we often suppress the norm $\lVert\cdot\rVert$ from the notation for a (semi)normed ring (respectively, Banach ring). For a seminormed abelian group $A$ we denote by $A_{\leq r}$ the closed ball 
(in the seminormed topology) of radius $r>0$ in $A$, i.e.,\begin{equation*}A_{\leq r}=\{\, f\in A\mid \lVert f\rVert\leq r\,\}.\end{equation*}Since we restrict our attention to nonarchimedean seminorms, this is always a subgroup of $A$, open with respect to the topology defined by the given seminorm, and if $A$ is a seminormed ring and $c\leq1$, then $A_{\leq c}$ is an open subring. 
We often abbreviate \(A^{\circ}\) for \(A_{\leq 1}\).
Moreover, we use the same notation for closed unit balls of (possibly higher-rank) valuations on a ring $A$, i.e., $A_{v\leq1}$ is the set $f\in A$ such that $v(f)\leq1$ for the valuation $v$ on $A$.

A seminorm on a ring $A$ is said to be power-multiplicative (resp., multiplicative) if $\lVert f^{n}
\rVert = \lVert f \rVert^{n}$ for all $f$ and all integers $n \geq 0$ (resp., $\lVert fg \rVert =
\lVert f \rVert \lVert g \rVert$ for all $f, g \in A$). A multiplicative norm on a ring $A$ is also
called an absolute value. By a valued field $K$ we mean a normed ring $(K, \lVert\cdot\lVert)$ whose underlying ring is a field and whose norm is an absolute value. If the valued field $K$ is moreover a Banach ring, we call it a nonarchimedean field. For a valued field $K$, we denote by $\widetilde{K}$ its residue field and by $\vert K^{\times}\vert$ its value group.  
\begin{defn}A valued field extension, or extension of valued fields, is a field extension $L/K$ between valued fields, with absolute values $v_{K}$, $v_{L}$, such that $v_{K}$ is the restriction of $v_{L}$ to $K$. In this case, we also call the valued field $L$ an extension of $K$. We call a valued field extension $L/K$ a nonarchimedean field extension if the valued fields $K$, $L$ are nonarchimedean fields. \end{defn}
We shall use the following conditions on nonarchimedean field extensions. 
\begin{defn}[Immediate extension]An extension of valued fields $K\subseteq L$ is called immediate if $L$ has the residue field and value group as $K$. A nonarchimedean field which admits no proper immediate extensions is called spherically complete.\end{defn}
\begin{defn}[Semi-immediate extension]An extension of valued fields $K\subseteq L$ is called a semi-immediate extension if the group $\vert L^{\times}\vert/\vert K^{\times}\vert$ is torsion (i.e., if $\vert L^{\times}\vert\subseteq \sqrt{\vert K^{\times}\vert}$) and the residue field of $L$ is algebraic over the residue field of $K$. \end{defn}
Note that a nonarchimedean field which is both algebraically closed and spherically complete admits
no proper semi-immediate extensions: Indeed, for any algebraically closed nonarchimedean field $K$,
the value group of $K$ is divisibly closed and the residue field of $K$ is algebraically closed.

For a nonarchimedean field $K$, the following class of nonarchimedean field extensions of $K$ plays
a prominent role in this paper. Following Ducros \cite{Ducros07}, 1.1, we call any $l$-tuple
$(r_1,\dots, r_l)\in\mathbb{R}_{>0}^{l}$ of positive real numbers a polyradius, and, given a
nonarchimedean field $K$, we call $(r_1,\dots, r_l)$ a $K$-free polyradius if the elements
$r_1,\dots, r_l$ define $\mathbb{Q}$-linearly independent elements of the ``multiplicative''
$\mathbb{Q}$-vector space $(\mathbb{R}_{>0}/\sqrt{\vert K^{\times}\vert})\otimes_{\mathbb{Z}}\mathbb{Q}$. We will say that the polyradius is free (instead of $K$-free) whenever the nonarchimedean field $K$ is understood from the context. 
\begin{defn}[The fields $K_{r_1,\dots, r_n}$]\label{Definition of Berkovich's fields}
    For a nonarchimedean field $K$ and a polyradius $(r_1,\dots, r_n)\in\mathbb{R}_{>0}^{n}$, we let \(K_{r_{1},\dots, r_{n}}\) be the completion of
    the rational function field \(K(x_{1},\dots, x_{n})\) with respect to the (canonical extension of the) 
    \((r_1,\dots, r_n)\)-Gauss norm \(\vert\cdot\vert_{r_{1},\dots, r_{n}}\).
    
    In the same vein, we let $K_{r_1,\dots, r_n}^{\perfd}$ be the completion of the field $K(x_{1}^{1/p^{\infty}},\dots, x_{n}^{1/p^{\infty}})$ with respect to the canonical extension of the $(r_1,\dots, r_n)$-Gauss norm on $K[x_{1}^{1/p^{\infty}}, \dots, x_{n}^{1/p^{\infty}}]$; note that if $K$ is a perfect(oid) of characteristic $p$, this is the completed perfection of $K_{r_1,\dots, r_n}$. \end{defn}
The nonarchimedean fields $K_{r_1,\dots, r_n}$ appear in Ducros' Nullstellnesatz for Berkovich-affinoid algebras (\cite{Ducros07}, Théorème 2.7) as well as in the description of topologically finitely generated nonarchimedean field extensions of $K$ (see Lemma~\ref{Temkin's remark and polyradii 3}), which plays an important philosophical role in this paper. For a single radius $r>0$, the field $K_{r}$ is the completed residue field of a type (II) point on the Berkovich affine line $\mathbb{A}_{K}^{1,\an}$ over $K$ if $r\in\sqrt{\vert K^{\times}\vert}$, and it is the completed residue field of a type (III) point if $r\not\in \sqrt{\vert K^{\times}\vert}$, cf.~\cite{Berkovich}, 1.4.4, or Proposition \ref{Classification of points} below.     
\begin{lemma}\label{Type 3 points}If the polyradius $(r_1,\dots, r_n)$ is free, then $K_{r_{1},\dots, r_{n}}$ and $K_{r_1,\dots, r_n}^{\perfd}$ are given by \[
    K_{r_{1},\dots, r_{n}} = \left\{ \sum_{i=(i_1,\dots, i_{n}) \in \mathbb{Z}^{n}}^{} a_{i}x^{i_{1}}\cdot~\dots~\cdot x^{i_{n}} \mid a_{i} \in K, a_{i}r^{i} \xrightarrow{} 0 \text{ as } \vert i\vert
    \xrightarrow{} \infty \right\}
    \] and \begin{equation*}K_{r_1,\dots, r_n}^{\perfd}=\left\{ \sum_{i=(i_1,\dots, i_n) \in \mathbb{Z}[1 / p]^{n}} a_{i}x_{1}^{i_{1}}\cdot~\dots~\cdot x_{n}^{i_{n}} \mid \text{ for every } \epsilon>0, \{i\mid \vert a_{i}\vert r_{1}^{i_{1}}\cdot~\dots~\cdot r_{n}^{i_{n}}\geq\epsilon\} \text{ is finite } \right\}.\end{equation*}\end{lemma}
\begin{proof}The proof is completely analogous to the proof of \cite{Gleason22N}, Lemma 4.6, but we include it here for the reader's convenience. However, we only spell out the proof for $K_{r_1,\dots, r_n}$; the argument for $K_{r_1,\dots, r_n}^{\perfd}$ is the same except that the index set is then $\mathbb{Z}[1/p]^{n}$ rather than $\mathbb{Z}^{n}$. 

Let $K\langle x_{1},\dots, x_{n}\rangle_{(r_1,\dots, r_n)}$ be the ring of convergent power series
on the right-hand side of the equality, equipped with its corresponding Gauss norm. We claim that
the freeness of the polyradius $(r_1,\dots, r_n)$ implies that this ring is a field. Indeed, if for
some \begin{equation*}f=\sum_{i\in\mathbb{Z}^{n}}a_{i}x_{1}^{i_{1}}\cdot~\dots~\cdot
    x_{n}^{i_{n}}\in K\langle x_1,\dots, x_n\rangle_{(r_1,\dots, r_n)}\end{equation*} the sequence
    \begin{equation*}\{\vert a_{i}\vert r_{1}^{i_{1}}\cdot~\dots~\cdot
        r_{n}^{i_{n}}\,\}_{i\in\mathbb{Z}^{n}}\end{equation*}attains its maximum at two indices $i,
        j\in\mathbb{Z}^{n}$, then $a_{j}\neq0$ and
        \begin{equation*}r_{1}^{i_{1}-j_{1}}\cdot\dots\cdot r_{n}^{i_{n}-j_{n}}\in \frac{\vert
            a_{j}\vert}{\vert a_{i}\vert}\in \sqrt{\vert K^{\times}\vert},\end{equation*}so by the assumption that the polyradius $(r_1,\dots, r_n)$ is free, we must have $i_{1}-j_{1}=0,\dots, i_{n}-j_{n}=0$, i.e., $i=j$. Thus the sequence $\{\, \vert a_{i}\vert r_{1}^{i_{1}}\cdot~\dots~\cdot r_{n}^{i_{n}}\,\}_{i\in\mathbb{Z}^{n}}$ attains its maximum at a unique index $N\in\mathbb{Z}^{n}$. Then for every $i\neq N$, the element \begin{equation*}(a_{N}^{-1}x_{1}^{-N_{1}}\cdot~\dots~\cdots x_{n}^{-N_{n}})(a_{i}x_{1}^{i_{1}}\cdot~\dots~x^{i_{n}})\end{equation*}has norm $<1$. It follows that $a_{N}^{-1}x_{1}^{-N_{1}}\cdot~\dots~\cdots x_{n}^{-N_{n}}f-1$ is topologically nilpotent and thus $f$ is a unit. 

Once we know that the ring $K\langle x_{1},\dots, x_{n}\rangle_{(r_1,\dots, r_n)}$ is a field, we know that it must contain the field of fractions $K(x_1,\dots, x_n)$ of $K[x_1,\dots, x_n]$, since it clearly contains $K[x_1,\dots, x_n]$. Since the norm on $K\langle x_{1},\dots, x_{n}\rangle_{(r_1,\dots, r_n)}$ is the $(r_1,\dots, r_n)$-Gauss norm, it must then also contain the completion $K_{r_1,\dots, r_n}$ of $K(x_1,\dots, x_n)$ with respect to the $(r_1,\dots, r_n)$-Gauss norm. But $K[x_1,\dots, x_n]$ (and, a fortiori, $K(x_1,\dots, x_n)$) is dense in $K\langle x_{1},\dots, x_{n}\rangle_{(r_1,\dots, r_n)}$, so $K\langle x_{1},\dots, x_{n}\rangle_{(r_1,\dots, r_n)}$ has to be equal to $K_{r_1,\dots, r_n}$. \end{proof}
\begin{rmk}By contrast, if all the radii $r_1,\dots, r_n$ belong to the divisible closure of the value group of $K$, then $K_{r_1,\dots, r_n}$ is topologically isomorphic (by rescaling) to $K_{1,\dots, 1}$, i.e., to the completion of the rational function field $K(x_1,\dots, x_n)$ with respect to the usual Gauss norm.

In particular, the above lemma need no longer be true if we omit the freeness condition on the polyradius. We note that in the original work of Berkovich \cite{Berkovich} and in Ducros \cite{Ducros07}, the notation $K_{r_1,\dots, r_n}$ always refers to the affinoid algebra (in the sense of Berkovich) given by the right hand side of the equality in Lemma \ref{Type 3 points}; this affinoid algebra may fail to be a field if the polyradius $(r_1,\dots, r_n)$ is not free. In this paper, we want to treat free and non-free polyradii on an equal footing and so we follow Temkin in \cite{Temkin04}, paragraph preceding Proposition 3.1, in using the notation $K_{r_1,\dots, r_n}$ as in Definition \ref{Definition of Berkovich's fields}, regardless of whether the polyradius $(r_1,\dots, r_n)$ is free or not. In particular, we employ the same notation $K_{r}$ both for completed residue fields of type (III) points on $\mathbb{A}_{K}^{1,\an}$ and for the completed residue fields of type (II) points.\end{rmk} 
\begin{lemma}For any nonarchimedean field $K$ and radius $r>0$, the field $K_{r}$ is not a semi-immediate extension of $K$.\end{lemma}
\begin{proof}If $r\in \sqrt{\vert K^{\times}\vert}$, i.e., $r^{k}=\vert c\vert$ for some $c\in
K^{\times}$ and some \(k \in \mathbb{N}\), then the residue field of $K_{r}$ has transcendence
degree $1$ over the residue field of $K$, with the image of $\frac{x^{k}}{c}$ as the new
transcendental variable. If $r\not\in \sqrt{\vert K^{\times}\vert}$, then the residue field of $K_{r}$ coincides with the residue field of $K$, but the group $\vert K_{r}^{\times}\vert/\vert K^{\times}\vert$ is free of rank $1$, with free generator $r$ (see also \cite{Kedlaya13}, Corollary 2.27(ii) and (iii)).\end{proof}
The following lemma allows us to relate the nonarchimedean fields $K_{r_1,\dots, r_n}$ to Abhyankar valuations on finitely generated field extensions of a nonarchimedean field $K$. 
\begin{lemma}\label{Temkin's remark and polyradii}Let $K$ be a nonarchimedean field and let $L/K$ be
    a valued field extension, with valuation $v$, such that $L$ is finitely generated and purely
    trascendental over $K$. Suppose that there exists a transcendence basis $B=\{x_1,\dots, x_n\}$
    of $L$ over $K$ with $v(x_{i})\leq1$ for all $i$ and some $k\in \{0, 1, \dots, n\}$ such that
    the images \begin{equation*}\widetilde{x_{1}},\dots, \widetilde{x_{k}}\in
        \widetilde{L}\end{equation*}of $x_1,\dots, x_{k}$ form a transcendence basis of
        $\widetilde{L}$ over $\widetilde{K}$ and the elements $v(x_{1}),\dots, v(x_{n})$ form a
        basis of the ``multiplicative'' $\mathbb{Q}$-vector space
        \begin{equation*}v(L^{\times})/v(K^{\times})\otimes_{\mathbb{Z}}\mathbb{Q}.\end{equation*}Then
    the completion $\widehat{L}$ of $L$ with respect to $v$ is of the form $K_{r_1,\dots, r_n}$,
with $r_{i}=v(x_{i})$ for all $i=1,\dots, n$.\end{lemma}
\begin{proof}We write $r_{i}=v(x_{i})$ for all $i$. Suppose first that we know the assertion to be true in the cases $k=0$ and $k=n$, and let $k\in\{1,\dots, n-1\}$. By the case $k=0$, the completion of $K(x_1,\dots, x_{k})$ (with respect to the restriction of $v$) is $K_{r_1,\dots, r_k}$. Then the completion $\widehat{L}$ of $L$ coincides with the completion $\widehat{L^{\prime}}$ of $L^{\prime}=K_{r_1,\dots, r_{k}}(x_{k+1},\dots, x_{n})$ with respect to the canonical extension of $v$ (which we again denote by $v$). By Abhyankar's inequality (\cite{BourbakiAlg}, Ch.~VI, \S10, \no 3, Cor.~1) for the finitely generated field extension $K(x_1,\dots, x_{k})/K$, we see that \begin{equation*}r_1,\dots, r_k\in \sqrt{v(K^{\times})}.\end{equation*}Hence the elements $v(x_{k+1}),\dots, v(x_{n})$ form a basis of the ``multiplicative'' $\mathbb{Q}$-vector space $v(L^{\prime\times})/\vert K_{r_1,\dots, r_{k}}^{\times}\vert\otimes_{\mathbb{Z}}\mathbb{Q}$. Applying the case $k=n$, to the valued field extension $L^{\prime}/K_{r_1,\dots, r_{k}}$, we see that the completion $\widehat{L}=\widehat{L^{\prime}}$ is given by $(K_{r_1,\dots, r_{k}})_{r_{k+1},\dots, r_{n}}=K_{r_1,\dots, r_n}$. 

It remains to treat the cases $k=0$ and $k=n$ separately. Note that proving the equality $\widehat{L}=K_{r_1,\dots, r_n}$ is equivalent to showing that the restriction of $v$ to the polynomial ring $K[x_1,\dots, x_n]$ is equal to the $(r_1,\dots, r_n)$-Gauss norm $\vert\cdot\vert_{r_1,\dots, r_n}$, i.e., that \begin{equation*}v(f)=\max_{i=(i_1,\dots, i_n)}\vert a_{i_{1},\dots, i_{n}}\vert r_{1}^{i_{1}}\cdot~\dots~\cdot r_{n}^{i_{n}}=\vert a_{i_{1},\dots, i_{n}}\vert v(x_{1})^{i_{1}}\cdot~\dots~\cdot v(x_{n})^{i_{n}}\end{equation*}for all $f\in K[x_1,\dots, x_n]$. In the case $k=0$, the elements $r_{1},\dots, r_{n}$ are linearly independent in the ``multiplicative'' $\mathbb{Q}$-vector space $v(L^{\times})/v(K^{\times})\otimes_{\mathbb{Z}}\mathbb{Q}$. By the same argument as in the second paragraph of the proof of Lemma~\ref{Type 3 points}, this implies that, for every element \begin{equation*}f=\sum_{i=(i_1,\dots, i_n)}a_{i}x_{1}^{i_{1}}\cdot~\dots~\cdot x_{n}^{i_{n}}\in K[x_1,\dots, x_n],\end{equation*}the set \begin{equation*}\{\, \vert a_{i}\vert r_{1}^{i_{1}}\cdot~\dots~\cdot r_{n}^{i_n}\mid i=(i_1,\dots, i_n)\in \mathbb{Z}_{\geq0}^{n}\,\}\end{equation*}has a unique maximum. By usual properties of nonarchimedean norms, this entails $v(f)=\max_{i}\vert a_{i}\vert r_{1}^{i_{1}}\cdot~\dots~\cdot r_{n}^{i_{n}}$, proving the lemma in the case $k=0$. 

In the case $k=n$, the assumptions $r_{i}=v(x_{i})\leq 1$ and $\widetilde{x_{i}}\neq0$ imply that $r_{i}=v(x_{i})=1$ for all $i=1,\dots, n$. We want to show that in this case the restriction of $v$ to the polynomial ring $K[x_1,\dots, x_n]$ must be equal to the usual Gauss norm $\lVert\cdot\rVert=\vert\cdot\vert_{1,\dots, 1}$. For every element $f\in K[x_1,\dots, x_n]$, we can find some $a\in K^{\times}$ and $m\in\mathbb{Z}_{>0}$ with $\lVert af^{m}\rVert=1$. We want to prove that $v(f)=\lVert f\rVert$. Since both $v$ and the Gauss norm $\lVert\cdot\rVert$ are multiplicative norms, we may replace $f$ with $af^{m}$, for $a\in K^{\times}$ and $m>0$ as above, and assume that $\lVert f\rVert=1$. We need to prove that $v(f)=1$ in this case. Write $f=\sum_{i=(i_1,\dots, i_n)}a_{i}x_{1}^{i_{1}}\cdot~\dots~\cdot x_{n}^{i_{n}}$. If $v(f)<1$, then we get the relation \begin{equation*}\sum_{i=(i_1,\dots, i_n)}\widetilde{a_{i}}\widetilde{x_{1}}^{i_{1}}\cdot~\dots~\cdot\widetilde{x_{n}}^{i_{n}}=0\end{equation*}in the residue field $\widetilde{L}$ of $L$. On other hand, at least one of the $\widetilde{a_{i}}$ must be non-zero by the assumption that $\lVert f\rVert=1$. This contradicts the assumption that \begin{equation*}\widetilde{x_{1}},\dots, \widetilde{x_{n}}\in \widetilde{L}\end{equation*}are algebraically independent, so $v(f)=1$. \end{proof}
\begin{lemma}\label{Temkin's remark and polyradii 2}Let $K$ be a nonarchimedean and let $L/K$ be a finitely generated valued field extension, with valuation $v$, and let \begin{equation*}n=\rk_{\mathbb{Q}}(v(L^{\times})/v(K^{\times}))+\trdeg(\widetilde{L}/\widetilde{K}).\end{equation*}There exists an algebraically independent subset $B=\{x_{1},\dots, x_{n}\}$ of $L/K$ of cardinality $n$ with $v(x_{i})\leq1$ for all $i=1,\dots, n$ such that, for some $k\in \{0, 1,\dots, n\}$, the images \begin{equation*}\widetilde{x_{1}},\dots, \widetilde{x_{k}}\in \widetilde{L}\end{equation*}of $x_1,\dots, x_{k}$ form a transcendence basis of $\widetilde{L}$ over $\widetilde{K}$ and the elements $v(x_{k+1}),\dots, v(x_{n})$ form a basis of the ``multiplicative'' $\mathbb{Q}$-vector space $v(L^{\times})/v(K^{\times})\otimes_{\mathbb{Z}}\mathbb{Q}$. In particular, by Lemma~\ref{Temkin's remark and polyradii}, the completion of $K(B)$ with respect to $v$ is of the form $K_{r_1,\dots, r_n}$. \end{lemma}
\begin{proof}Let $k$ be the transcendence degree of the residue field extension $\widetilde{L}/\widetilde{K}$, so \begin{equation*}\rk_{\mathbb{Q}}(v(L^{\times})/v(K^{\times}))=n-k.\end{equation*}Choose a transcendence basis $\widetilde{x_{1}},\dots, \widetilde{x_{k}}$ of $\widetilde{L}$ over $\widetilde{K}$ and lift it to elements $x_1,\dots, x_k\in L^{\circ}$ (note that we necessarily have $v(x_{i})=1$ for all $i=1,\dots, k$). The elements $x_1,\dots, x_k$ are algebraically independent in $L$. Indeed, suppose that \begin{equation*}\sum_{i=(i_1,\dots, i_{k})}a_{i}x_{1}^{i_{1}}\cdot~\dots~\cdot x_{k}^{i_{k}}=0\end{equation*}is an algebraic relation with coefficients in $K$. Choose $i$ such that $\vert a_{i}\vert$ is maximal. Multiplying by above equation by $a_{i}^{-1}$, we obtain an algebraic relation all whose coefficients have absolute value $\leq1$ and at least one of whose coefficients has absolute value $1$. Applying the reduction map $L^{\circ}\to \widetilde{L}$, we get a non-trivial algebraic relation between $\widetilde{x_1},\dots, \widetilde{x_{k}}$, contradicting the algebraic independence of $\widetilde{x_{1}}, \dots, \widetilde{x_{k}}$ over $\widetilde{K}$. 

Let $x_{k+1},\dots, x_{n}$ be elements of $L$ such that $r_{k+1}=v(x_{k+1}),\dots, r_{n}=v(x_{n})$ form a basis of the ``multiplicative'' $\mathbb{Q}$-vector space $v(L^{\times})/v(K^{\times})\otimes_{\mathbb{Z}}\mathbb{Q}$. The elements $x_{k+1},\dots, x_{n}$ are algebraically independent. Indeed, if they were not, then there would exist an algebraic relation \begin{equation*}\sum_{i=(i_{k+1},\dots, i_{n})}a_{i}x_{k+1}^{i_{k+1}}\cdot~\dots~\cdot x_{n}^{i_{n}}=0.\end{equation*}But since $r_{k+1},\dots, r_{n}$ are linearly independent in $v(L^{\times})/v(K^{\times})\otimes_{\mathbb{Z}}\mathbb{Q}$ (see the argument in the proof of Lemma~\ref{Type 3 points}), the set \begin{equation*}\{\vert a_{i}\vert r_{k+1}^{i_{k+1}}\cdot~\dots~\cdot r_{n}^{i_{n}}\mid i=(i_{k+1},\dots, i_{n})\in\mathbb{Z}^{n-k}\}\end{equation*}has a unique maximum, say $\vert a_{j}\vert r_{k+1}^{j_{k+1}}\cdot~\dots~\cdot r_{n}^{j_{n}}$, and the above relation implies $\vert a_{j}\vert=0$. It follows that $\vert a_{i}\vert=0$ for all $i$, that is, $a_{i}=0$ for all $i$. 

Finally, we show that the whole set $\{x_{1},\dots, x_{n}\}$ is algebraically independent. By way of contradiction, suppose that it is algebraically dependent. Since have seen that each of the sets $\{x_1,\dots, x_k\}$ and $\{x_{k+1},\dots, x_{n}\}$ is algebraically independent, this entails the existence of an algebraic relation \begin{equation*}\sum_{i=(i_1,\dots, i_{k})}a_{i}x_{1}^{i_{1}}\cdot~\dots~\cdot x_{k}^{i_{k}}=\sum_{j=(j_{k+1},\dots, j_{n})}b_{j}x_{k+1}^{j_{k+1}}\cdot~\dots~\cdot x_{n}^{j_{n}}\end{equation*}where neither all of the $a_{i}$ nor all of the $b_{j}$ are zero. By Lemma~\ref{Temkin's remark and polyradii}, the restriction of $v$ to $K(x_1,\dots, x_k)$ is the usual Gauss norm (for the variables $x_1,\dots, x_{k}$), so there exists an index $i^{\prime}$ such that \begin{equation*}v(\sum_{i}a_{i}x_{1}^{i_{1}}\cdot~\dots~\cdot x_{k}^{i_{k}})=\vert a_{i^{\prime}}\vert.\end{equation*}Using once again the linear independence of $r_{k+1},\dots, r_{n}$ in $v(L^{\times})/v(K^{\times})\otimes_{\mathbb{Z}}\mathbb{Q}$, the set \begin{equation*}\{\vert b_{j}\vert r_{k+1}^{j_{k+1}}\cdot~\dots~\cdot r_{n}^{j_{n}}\mid j=(j_{k+1},\dots, j_{n})\in \mathbb{Z}^{n-k}\}\end{equation*}has a unique maximum, and then, by the usual properties of nonarchimedean norms, this unique maximum is equal to \begin{equation*}v(\sum_{j=(j_{k+1},\dots, j_{n})}b_{j}x_{k+1}^{j_{k+1}}\cdot~\dots~\cdot x_{n}^{j_{n}})=v(\sum_{i=(i_1,\dots, i_{k})}a_{i}x_{1}^{i_{1}}\cdot~\dots~\cdot x_{k}^{i_{k}})=\vert a_{i^{\prime}}\vert.\end{equation*}This implies that one of the $r_{k+1}, \dots, r_{n}$ belongs to the divisible closure of the value group $v(K^{\times})$, a contradiction. Thus the elements $x_1,\dots, x_n$ are algebraically independent over $K$, as required.  \end{proof}
We conclude this section with a variant of \cite{Temkin13}, Remark~2.1.2, which will play a crucial role in the proof of our main result. We call a finitely generated valued field extension $L/K$, with valuation $v$, Abhyankar if Abhyankar's inequality for this valued field extension is an equality, i.e., \begin{equation*}\rk_{\mathbb{Q}}(v(L^{\times})/v(K^{\times}))+\trdeg(\widetilde{L}/\widetilde{K})=\trdeg(L/K).\end{equation*}In this case, we also say that the valuation $v$ on $L$ is Abhyankar (over $K$). Note that our semi-immediate valued field extensions are called transcendentally immediate in \cite{Temkin13}.   
\begin{lemma}\label{Temkin's remark and polyradii 3}Let $K$ be a nonarchimedean field and let $L/K$ be a finitely generated valued field extension, with valuation $v$, and let \begin{equation*}n=\rk_{\mathbb{Q}}(v(L^{\times})/v(K^{\times}))+\trdeg(\widetilde{L}/\widetilde{K}).\end{equation*}There exists an algebraically independent subset $\{x_1,\dots, x_n\}$ of $L/K$ of cardinality $n$ with $v(x_{i})\leq 1$ such that the restriction of $v$ to $K(x_1,\dots, x_n)$ is Abhyankar and the extension $L/K(x_1,\dots, x_n)$ is semi-immediate.

In particular, there exists a polyradius $(r_1,\dots, r_n)\in (0, 1]^{n}$ such that the nonarchimedean field extension $\widehat{L}/K$ factors as \begin{equation*}\widehat{L}/K_{r_1,\dots, r_n}/K\end{equation*}and the nonarchimedean field extension $\widehat{L}/K_{r_1,\dots, r_n}$ is semi-immediate.\end{lemma}
\begin{proof}Let $x_1,\dots, x_n$ be chosen as in Lemma~\ref{Temkin's remark and polyradii 2}. Then the finitely generated valued field extension $K(x_1,\dots, x_n)/K$ is Abhyankar and the extension $L/K(x_1,\dots, x_n)$ is semi-immediate by fiat. By Lemma~\ref{Temkin's remark and polyradii}, the completion of $K(x_1,\dots, x_n)$ with respect to the restriction of $v$ is of the nonarchimedean field $K_{r_1,\dots, r_n}$, where $r_{i}=v(x_{i})$ for all $i=1,\dots, n$. Since the valued field extension $L/K(x_1,\dots, x_n)$ is semi-immediate, so is the nonarchimedean field extension $\widehat{L}/K_{r_1,\dots, r_n}$ obtained by taking completions. \end{proof}
\begin{cor}\label{Temkin's remark 4}Let $L/K$ be a perfectoid field extension of characteristic $p$ which is perfectly topologically finitely generated, i.e., $L$ has a dense subfield which is the perfection of a finitely generated field extension of $K$. Then $L/K$ factors as \begin{equation*}L/K_{r_1,\dots, r_l}^{\perfd}/K\end{equation*}for a polyradius $(r_1,\dots, r_l)$ and a semi-immediate extension $L/K_{r_1,\dots, r_l}^{\perfd}$.\end{cor}
\begin{proof}Let $L'$ be a dense subextension of $L/K$ which is the perfection of a finitely generated extension $L''/K$. By Lemma~\ref{Temkin's remark and polyradii 3}, there exists an algebraically independent subset $B$ of $L''$ such that the completion of $K(B)$ is of the form $K_{r_1,\dots, r_l}$, $l=\vert B\vert$, and such that $L''/K(B)$ is a semi-immediate extension. But then $L/K_{r_1,\dots, r_l}^{\perfd}$ is also a semi-immediate extension since completed perfections are semi-immediate.\end{proof}  

\subsection{Berkovich's classification of points}

Let $K$ be a nonarchimedean field. For $a \in \widehat{\overline{K}}$ in the completed algebraic closure $\widehat{\overline{K}}$ of $K$ we denote by $B(a, r)$ the closed disk of radius $r$ around $a$, i.e., \begin{equation*}B(a, r)=\{\, b \in \widehat{\overline{K}} \mid \vert b - a\vert \leq r \,\}.\end{equation*}
Consider the polynomial ring $K[T]$. For every closed disk $B(a, r) \subset \widehat{\overline{K}}$ as above we define a multiplicative seminorm $\vert\cdot\vert_{B(a, r)}$ on $K[T]$ by 
\begin{equation*} \vert f\vert_{B(a,r)} := \sup_{z \in B(a, r)} \vert f(z)\vert \ \ \ (f \in K[T]).\end{equation*}Note that this seminorm extends the given absolute value on $K$. We recall Berkovich's classification of points on the Berkovich affine line $\mathbb{A}_{K}^{1,\an}$ (\cite{Berkovich}, 1.4.4), i.e., of multiplicative seminorms on $K[T]$ which extend the given absolute value on the nonarchimedean field $K$. For a point $x$ on a Berkovich-analytic space $X$ over a nonarchimedean field $K$ (such as $\mathbb{A}_{K}^{1,\an}$ or such as the Gelfand spectrum $\mathscr{M}(A)$ of an affinoid algebra $A$ in the sense of Berkovich), we denote by $\mathcal{H}(x)$ the completed residue field of $X$ at $x$; for an affinoid neighborhood $\mathscr{M}(A)$ of $x$, the point $x$ is a continuous rank-$1$ valuation on $A$ (over $K$) and $\mathcal{H}(x)$ is the completion of the field of fractions $Q(A/\ker(x))$ of $A/\ker(x)$ with respect to the absolute value induced by $x$. 
\begin{prop}[Berkovich \cite{Berkovich}, 1.4.4, see also \cite{Kedlaya13}, Corollary 2.27]\label{Classification of points}Every point $x \in\mathbb{A}_{K}^{1,\an}$ belongs to one of the following categories:\begin{enumerate}[(1)]
\item Points of type $(I)$: There exists a point $a \in \widehat{\overline{K}}$ such that, for all $f \in K[T]$ we have $\vert f\vert_{x}=\vert f(a)\vert$. The completed residue field $\mathcal{H}(x)$ is a finite extension of $K$. In particular, it is a semi-immediate extension of $K$. 
\item Points of type $(II)$: There exist a point $a \in \widehat{\overline{K}}$ and a closed disk
    $B(a, r) \subset \widehat{\overline{K}}$, with radius $r$ in $\sqrt{\vert K^{\times}\vert}$, such that \begin{equation*}\vert f\vert_{x}=\vert f\vert_{B(a,r)}\end{equation*}for all $f \in K[T]$. The completed residue field is of the form $\mathcal{H}(x)=K_{r}$ with $r\in \sqrt{\vert K^{\times}\vert}$. In this case, $\vert \mathcal{H}(x)^{\times}\vert/\vert K^{\times}\vert$ is a finite group and $\widetilde{\mathcal{H}(x)}$ is finitely generated over $\widetilde{K}$ of transcendence degree $1$.
\item Points of type $(III)$: The same as for points of type (II), but the radius $r$ of $B(a, r)$
    does not belong to $\sqrt{\vert K^{\times}\vert}$. The completed residue field is $\mathcal{H}(x)=K_{r}$ with $r\not\in \sqrt{\vert K^{\times}\vert}$. The group $\vert \mathcal{H}(x)^{\times}\vert/\vert K^{\times}\vert$ is a finitely generated abelian group of rank $1$, generated by $r$.    
\item Points of type $(IV)$: There exists a nested sequence of closed disks \begin{equation*}B(a_1, r_1) \supseteq B(a_2, r_2) \supseteq B(a_3, r_3) \supseteq \dots \end{equation*}in $\widehat{\overline{K}}$ with empty intersection such that \begin{equation*} \vert f\vert_{x} = \lim_{n}\vert f\vert_{B(a_n,r_n)}=\inf_{n} \vert f\vert_{B(a_n,r_n)} \ (f \in K[T]) 
\end{equation*}In this case, the completed residue field $\mathcal{H}(x)$ is a semi-immediate extension of $K$. 
\end{enumerate}If the completed algebraic closure $\widehat{\overline{K}}$ is spherically complete, then only points of types (I), (II) and (III) occur.\end{prop} 

Moreover, we have the same description for the points in the Berkovich closed unit ball, i.e.
\(\mathscr{M}(K\langle T\rangle )\), and, more generally, in a closed ball $\mathscr{M}(K\langle r^{-1}T\rangle)$ of radius $r>0$, subject to the condition that the ball only contains
points which are dominated by the $r$-Gauss norm.

\subsection{Shilov Boundaries}
\label{sec:shilov:boundary}

In this subsection, we fix a nonarchimedean field $K$. For a seminormed ring $(A, \lVert\cdot\rVert)$, we consider the spectral seminorm \begin{equation*}\vert f\vert_{\spc}=\inf_{n}\lVert f^{n}\rVert^{1/n}=\lim_{n\to\infty}\lVert f^{n}\rVert^{1/n}, f\in A.\end{equation*}If $A$ is a seminormed algebra over $K$, then this definition does not depend on the choice of the $K$-algebra seminorm $\lVert\cdot\rVert$ defining the topology on $A$: Indeed, any two such $K$-algebra seminorms $\lVert\cdot\rVert$, $\lVert\cdot\rVert'$ have to be bounded-equivalent, i.e., \begin{equation*}\lVert \cdot\rVert C\leq \lVert\cdot\rVert'\leq C'\lVert\cdot\rVert\end{equation*}for some constants $C, C'>0$, and then their spectral seminorms, defined as above, have to be equal by \cite{BGR}, Proposition 1.3.1/2. We note that $\vert\cdot\vert_{\spc}$ is always the maximal power-multiplicative seminorm on $A$ bounded by the original seminorm $\lVert\cdot\rVert$.

Let us also recall the definition of the reduction of a seminormed ring.
\begin{defn}[Reduction of a seminormed ring]For any seminormed ring $A$, we define the reduction $\widetilde{A}$ of $A$ to be \begin{equation*}\widetilde{A}=A^{\circ}/A^{\circ\circ}.\end{equation*}\end{defn}
If the seminorm on $A$ is power-multiplicative, then this is the quotient of the closed unit ball by the open unit ball. Note that this notation is compatible with our notation for residue fields of nonarchimedean fields.


Recall from \cite{Berkovich}, Theorem 1.3.1, that the
spectral seminorm can also be described as \begin{equation*}\vert
f\vert_{\spc}=\sup_{v\in\mathscr{M}(A)}v(f)=\max_{v\in\mathscr{M}(A)}v(f),\end{equation*}where the
supremum is a maximum since the Gelfand spectrum $\mathscr{M}(A)$ is a compact Hausdorff space (\cite{Berkovich}, Theorem 1.2.1). It is then reasonable to ask for minimal closed subsets $\mathcal{S}$ of $\mathscr{M}(A)$ on which every element $f\in A$ achieves is maximum $\vert f\vert_{\spc}$.      
\begin{defn}[Shilov boundary: \cite{Guennebaud}, \S4, \cite{Dine25}, Definition 2.25]For a seminormed $K$-algebra $A$, a subset $\mathcal{S}\subset\mathscr{M}(A)$ is called a boundary for $A$ if every element $f\in A$ attains its maximum $\vert f\vert_{\spc}$ on $\mathcal{S}$. If there exists a unique minimal closed boundary $\mathcal{S}$, it is called the Shilov boundary for $A$.\end{defn}
It is known (see \cite{Guennebaud}, Ch.~I, Théorème 3, \cite{EM1}, Theorem C) that the Shilov boundary exists for any seminormed $K$-algebra $A$ (with suitable definitions, it actually exists for any seminormed ring and even any seminormed monoid, see the work of Guennebaud \cite{Guennebaud}). In this paper, we only need the notion of Shilov boundary for affinoid algebras in the sense of Tate. In this case, the Shilov boundary has an interesting description, which due to Berkovich in \cite{Berkovich}. To recall this description, we introduce
some notation. For any affinoid algebra in the sense of Berkovich $A$, we consider the map
\begin{equation*}\spc: \Spa(A, A^{\circ})\to
    \Spf(A^{\circ})=\Spec(\widetilde{A})\end{equation*}which sends a continuous valuation on $A$ to
    its center $\spc(v)$ on $A^{\circ}$, i.e., \begin{equation*}\spc(v)=\{\, f\in A^{\circ}\mid
    v(f)<1\,\}.\end{equation*}This map is continuous; by contrast, its restriction to the Gelfand spectrum $\mathscr{M}(A)$ is anti-continuous for the usual compact Hausdorff topology on $\mathscr{M}(A)$ (the pre-image of a Zariski-closed subset of $\Spf(A^{\circ})$ is open in $\mathscr{M}(A)$). For any ring $R$, we denote by $\Min(R)$ the set of minimal prime ideals of $R$. We can now state Berkovich's description of the Shilov boundary for strictly affinoid algebras.
\begin{prop}[\cite{Berkovich}, Proposition 2.4.4(ii), (iii), cf.~\cite{Bhatt-Hansen}, Proposition 2.2 for the formulation using the adic spectrum]\label{Berkovich's result}
Let $A$ be a strictly $K$-affinoid algebra in the sense of Berkovich (i.e, an affinoid algebra in the sense of Tate). Then the pre-image $\spc^{-1}(\{\mathfrak{p}\})$ of any $\mathfrak{p}\in \Min(\widetilde{A})$ is a single rank-$1$ point and the set of such pre-images \begin{equation*}\spc^{-1}(\Min(\widetilde{A}))\end{equation*}is the Shilov boundary for $A$. In particular, the Shilov boundary is finite.\end{prop}
The Shilov boundary of a general Berkovich-affinoid algebra $A$ need not consist of pre-images of minimal prime ideals of the reduction as in the above proposition. However, Temkin established in \cite{Temkin04} an analog of Berkovich's description for all affinoid algebras using his notation of graded reduction, see \cite{Temkin04}, Proposition 3.3.

\section{Abhyankar valuations}
\label{sec:abhyankar:val}

Recall that any rank-$1$ valuation $v$ on a finitely generated field extension $L/K$ satisfies Abhyankar's inequality \begin{equation*}\rk_{\mathbb{Q}}(v(L^{\times})/v(K^{\times}))+\trdeg(\widetilde{L}/\widetilde{K})\leq \trdeg(L/K),\end{equation*}where the residue fields are the residue fields of $L$ and $K$ with respect to $v$ and its restriction $v\vert_{K}$ (\cite{BourbakiAlg}, Ch.~VI, \S10, \no 3, Cor.~1). Recall also (for example, from \cite{Knaf-Kuhlmann05}, \cite{Cutkosky22}) that the valuation $v$ is called Abhyankar if the above inequality is an equality. In this case, we also say that the finitely generated valued field extension $L/K$ is Abhyankar.    

Besides the notion of Abhyankar valuations on finitely generated field extensions, we also study a
related notion of (continuous) Abhyankar valuations on affinoid algebras (in the sense of Berkovich), which corresponds
precisely to the notion of Abhyankar points on Berkovich analytic spaces introduced by Ducros in
\cite{Ducros18}, 1.4.10. 
Following \cite{Ducros18}, 1.2.4 and A.4.11, we define, for any valued field extension $L/K$, \begin{equation*}d_{K}(L) \colonequals
\rk_{\mathbb{Q}}(v(L^{\times})/v(K^{\times}))+\trdeg(\widetilde{L}/\widetilde{K}).\end{equation*}
Furthermore, as in \cite{Ducros18}, 1.4.6, we use the notation \begin{equation*}d_{K}(v)\colonequals d_{K}(\mathcal{H}(v))\end{equation*}for a point $v$ on a Berkovich space $X$, where, as usual, $\mathcal{H}(v)$ denotes the completed residue field of $X$ at $v$. Recall that if $X=\mathscr{M}(A)$ for some Berkovich-affinoid algebra $A$, then $\mathcal{H}(v)$ is the completion of the field of fractions $Q(A/\ker(v))$ of $A/\ker(v)$ with respect to the absolute value induced by the bounded multiplicative seminorm $v: A\to\mathbb{R}_{\geq0}$. 
Note that a bounded multiplicative seminorm on $A$ is the same thing as a continuous rank-$1$ valuation on $A$ over $K$, i.e., a continuous rank-$1$ valuation which restricts to the given absolute value on the ground field $K$.  
\begin{defn}[Abhyankar valuation on an affinoid algebra]\label{Abhyankar valuations, definition}Let
    $K$ be a nonarchimedean field and let $A$ be a $K$-affinoid algebra in the sense of Berkovich.
    We call a continuous rank-$1$ valuation on $A$ over $K$ an Abhyankar valuation (over $K$) if
    \begin{equation*}d_{K}(v)=\inf_{V}\dim(A_{V}),\end{equation*}where $V$ ranges over the affinoid domains of $\mathscr{M}(A)$ containing $v$ and where $A_{V}$ is the affinoid algebra corresponding to an affinoid domain $V$.\end{defn}
By a result of Ducros, this definition simplifies considerably if $\mathscr{M}(A)$ is purely $d$-dimensional, i.e., if every irreducible component of $\mathscr{M}(A)$ has the same dimension $d$, in which case $d$ is equal to the dimension of $A$. 
\begin{lemma}\label{Abhyankar valuations, simplified}Let $A$ be a $K$-affinoid algebra in the sense
    of Berkovich over a nonarchimedean field $K$ and suppose that its Gelfand spectrum
    $\mathscr{M}(A)$ is purely $d$-dimensional. Then, for every affinoid domain $V$ in
    $\mathscr{M}(A)$, the corresponding affinoid algebra $A_{V}$ satisfies $\dim(A_{V})=\dim(A)=d$.
    In particular, a continuous rank-$1$ valuation on $A$ over $K$ is Abhyankar if and only if
    \begin{equation*}d_{K}(v)=\dim(A).\end{equation*}\end{lemma}
\begin{proof}By \cite{Ducros07}, Lemme 1.12, the dimension of $A_{V}$, for any affinoid domain $V$, is equal to the maximum of the dimensions of those irreducible components of $\mathscr{M}(A)$ which have non-empty intersection with $V$. Since $\mathscr{M}(A)$ is assumed to be purely $d$-dimensional, we see that $\dim(A_{V})=\dim(A)=d$. \end{proof}
\begin{cor}\label{Abhyankar valuations, simplified 2}Let $A$ be a $K$-affinoid algebra in the sense
    of Berkovich over a nonarchimedean field $K$ and suppose that $A$, as a ring without topology,
    is an integral domain. Then, for every affinoid domain $V$ of $\mathscr{M}(A)$, we have
    $\dim(A_{V})=\dim(A)$. In particular, a continuous rank-$1$ valuation $v$ on $A$ (over $K$) is
    Abhyankar if and only if \begin{equation*}d_{K}(v)=\dim(A).\end{equation*}\end{cor}
\begin{proof}If $A$ is an integral domain, $\mathscr{M}(A)$ has only one irreducible component, so the assertion follows from Lemma \ref{Abhyankar valuations, simplified}. \end{proof}
\begin{cor}\label{Abhyankar absolute values}Let $A$ be a $K$-affinoid algebra in the sense of
    Berkovich over a nonarchimedean field $K$ and let $v$ be a continuous absolute value on $A$ over
    $K$. Then $v$ is Abhyankar if and only if \begin{equation*}d_{K}(v)=\dim(A).\end{equation*}\end{cor}
\begin{proof}Any ring on which there exists an absolute value is an integral domain. \end{proof}
For strictly affinoid algebras which are not integral domains, we still have the following property.
\begin{lemma}\label{The kernel of an Abhyankar valuation}Let $A$ be a strictly $K$-affinoid algebra over a nonarchimedean field $K$ and let $v$ be a continuous rank-$1$
    valuation on $A$ over $K$. If $v$ is Abhyankar, then
    \begin{equation*}\height(\ker(v))+d_{K}(v)=\dim(A).\end{equation*}\end{lemma}
\begin{proof}By \cite{Ducros09}, Théorème 2.13 (in the strictly affinoid case, this theorem is due to Kiehl,
    \cite{Kiehl69}, Theorem 3.3), $A$ is excellent and, in particular, catenary, so if
    $\height(\ker(v))>0$, then $\dim(A/\ker(v))<\dim(A)$. By Abhyankar's inequality
    (\cite{BourbakiAlg}, Ch.~VI, \S10, \no 3, Cor.~1), \begin{equation*}d_{K}(v)\leq \inf_{L/K}\trdeg(L/K),\end{equation*}where $L/K$ ranges over the dense finitely generated subextensions of $Q(A/\ker(v))/K$, and, by Noether normalization for the strictly affinoid algebra $A/\ker(v)$, the right hand side of this inequality is $\leq\dim(A/\ker(v))$. It follows that $d_{K}(v)<\dim(A)$ if $\ker(v)$ is not a minimal prime ideal of $A$.\end{proof}

\subsection{Rational Abhyankar and weakly Shilov valuations}
        
Lemma~\ref{Temkin's remark and polyradii 3} shows that Abhyankar valuations generalize both
points of type (II) and points of type (III) in Berkovich's classification of points on the
nonarchimedean-analytic affine line (both have a residue field of the form described by
Lemma~\ref{Temkin's remark and polyradii 3}). This suggests a differentiation of Abhyankar valuations into two types, which generalize type (II) points and type (III) points, respectively.
\begin{defn}[Rational and irrational Abhyankar valuations]\label{Rational and irrational}We call an
    Abhyankar valuation $v$ on a Berkovich-affinoid algebra $A$ \textit{rational} if the group $\vert \mathcal{H}(v)^{\times}\vert/\vert K^{\times}\vert$ is torsion, i.e., if \begin{equation*}\vert \mathcal{H}(v)^{\times}\vert\subseteq \sqrt{\vert K^{\times}\vert}.\end{equation*}Equivalently, the Abhyankar valuation $v$ on $A$ is rational if and only if \begin{equation*}\trdeg(\widetilde{\mathcal{H}(v)}/\widetilde{K})=\inf_{V}\dim_{V}A_{V},\end{equation*}where $V$ ranges over the affinoid domains of $\mathscr{M}(A)$ containing $v$ and $A_{V}$ denotes the affinoid algebra of the affinoid domain $V$. 

We call the Abhyankar valuation $v$ \textit{irrational} if it is not rational.\end{defn}
\begin{lemma}\label{Temkin's remark and polyradii 4}Let $v$ be a continuous rank-$1$ valuation on $K\langle T_1,\dots, T_n\rangle$ over $K$. If $v$ is Abhyankar (respectively, rational Abhyankar, respectively, irrational Abhyankar), then \begin{equation*}\mathcal{H}(v)=K_{r_1,\dots, r_n}\end{equation*}for some polyradius $(r_1,\dots, r_n)\in (0, 1]^{n}$ (respectively, for some polyradius $(r_1,\dots, r_n)\in (0, 1]^{n}$ with $r_1,\dots, r_n\in \sqrt{\vert K^{\times}\vert}$, respectively, for some polyradius $(r_1,\dots, r_n)\in (0, 1]^{n}$ with $r_{i}\not\in \sqrt{\vert K^{\times}\vert}$ for at least one index $i\in \{1,\dots, n\}$). \end{lemma}
\begin{proof}By abuse of notation, we again denote by $v$ the rank-$1$ valuation on $\mathcal{H}(v)$. Note that $\mathcal{H}(v)$ is the completion, with respect to the valuation defined by $v$, of the rational function field $K(T_1,\dots, T_n)$. Since $v$ is Abhyankar, the valued field extension $K(T_1,\dots, T_n)/K$ is Abhyankar in the sense that Abhyankar's inequality for it is an equality. Applying Lemma~\ref{Temkin's remark and polyradii 3} to the purely transcendental valued extension $K(T_1,\dots, T_n)/K$, we see that there exists a transcendence basis $\{x_1,\dots, x_n\}$ of $K(T_1,\dots, T_n)/K$ such that the completion $\mathcal{H}(v)$ of $K(T_1,\dots, T_n)=K(x_1,\dots, x_n)$ with respect to $v$ is of the form $K_{r_1,\dots, r_n}$, where $r_{i}=v(x_{i})\leq1$ for all $i=1,\dots, n$. Moreover, if $v$ is rational Abhyankar (respectively, irrational Abhyankar), then $\vert \mathcal{H}(v)^{\times}\vert\subseteq \sqrt{\vert K^{\times}\vert}$ (respectively, $\vert \mathcal{H}(v)^{\times}\vert\not\subseteq \sqrt{\vert K^{\times}\vert}$) and thus $r_1,\dots, r_n\in \sqrt{\vert K^{\times}\vert}$ (respectively, there exists some $i\in\{1,\dots, n\}$ with $r_{i}\not\in \sqrt{\vert K^{\times}\vert}$).  \end{proof} 
In the case of strictly affinoid algebras, the notion of rational Abhyankar valuations turns out to be equivalent to the following notion which was first introduced by Bhatt and Hansen in \cite{Bhatt-Hansen}.
\begin{defn}[Weakly Shilov valuation, see \cite{Bhatt-Hansen}, Definition 2.5]Let $A$ be an affinoid algebra in the sense of Tate over a nonarchimedean field $K$. A continuous rank-$1$ valuation $v$ on $A$ over $K$ is called Shilov if it is an element of the Shilov boundary for $A$, and it is called weakly Shilov if there exists an open affinoid subspace $\Spa(B, B^{+})$ of $\Spa(A, A^{\circ})$ containing $v$ such that $v$ belongs to the Shilov boundary for $B$.\end{defn}
Using a result of Berkovich on the behaviour of the Shilov boundary under passage to affinoid subdomains, \cite{Berkovich}, Proposition 2.5.20, we can describe the above definition in multiple equivalent ways.  
\begin{lemma}\label{Weakly Shilov implies rationally weakly Shilov}Let $K$ be a nonarchimedean field and let $A$ be an affinoid algebra in the sense of Tate over $K$. For a continuous rank-$1$ valuation $v$ over $K$ the following are equivalent: \begin{enumerate}[(1)]\item $v$ is weakly Shilov. \item There exists a rational subset $\Spa(C, C^{+})$ of $\Spa(A, A^{\circ})$ with $v\in \Spa(C, C^{+})$ such that $v$ belongs to the Shilov boundary for $C$. \item There exists a strictly affinoid neighborhood $V$ of $v$ in $\mathscr{M}(A)$ such that $v$ belongs to the Shilov boundary for $A_{V}$. \item There exists an affinoid neighborhood $V$ of $v$ in $\mathscr{M}(A)$ such that $v$ belongs to the Shilov boundary for $A_{V}$.\end{enumerate}\end{lemma}
\begin{proof}The implications (2)$\Rightarrow$ (3) $\Rightarrow$ (4) and (2) $\Rightarrow$ (1) are obvious. We prove the implications (1)$\Rightarrow$(2) and (4)$\Rightarrow$ (2). 

(1) $\Rightarrow$ (2): Let $\Spa(B, B^{+})$ be an affinoid open neighborhood of $v$ in $\Spa(A, A^{\circ})$ such that $v$ belongs to the Shilov boundary for $B$. Let $\varphi$ denote the map $A\to B$. Since rational subsets form a basis for the topology of $\Spa(A, A^{\circ})$, there exists a rational subset $\Spa(C, C^{+})$ of $\Spa(A, A^{\circ})$ containing $v$ and contained in $\Spa(B, B^{+})$. Then $\Spa(C, C^{+})$ is also a rational subset of $\Spa(B, B^{+})$. Indeed, the map $(A, A^{\circ})\to (C, C^{+})$ factors through $\varphi$, so $(C, C^{+})$ is the rational localization of $(A, A^{\circ})$ given by elements $f_1,\dots, f_n, g\in A$ generating the unit ideal of $A$, then it is also the rational localization of $(B, B^{+})$ given by their images $\varphi(f_1),\dots, \varphi(f_n), \varphi(g)$. Since $\Spa(B, B^{+})$ is an affinoid open subspace of $\Spa(A, A^{\circ})$, the Tate ring $B$ is topologically of finite type over $K$, i.e., it is an affinoid algebra in the sense of Tate. Then the rational localization $B\to C$ corresponds to a (strictly) affinoid domain $V$ of $\mathscr{M}(B)$ containing $v$, with affinoid algebra $B_{V}=C$. Since $v$ belongs to the Shilov boundary for $B$, it then follows from \cite{Berkovich}, Proposition 2.5.20, that $v$ also belongs to the Shilov boundary for $C$. Thus $\Spa(C, C^{+})$ is the sought-for rational subset of $\Spa(A, A^{\circ})$.

(4)$\Rightarrow$ (2): Let $V$ be an affinoid neighborhood of $v$ such that $v$ belongs to the Shilov boundary for $A_{V}$. Since $A$ is strictly affinoid, the strictly affinoid neighborhoods of $v$ form a fundamental system of closed neighborhoods of $v$, by \cite{Berkovich}, Proposition~2.2.3. Hence we can find a strictly affinoid neighborhood $V'$ of $v$ contained in $V$. By the Gerritzen-Grauert Theorem, \cite{BGR}, Corollary~7.3.5/3, $V'$ is a union of finitely many strictly rational subdomains. In particular, there exists a strictly rational neighborhood $U$ of $v$ contained in $V'\subseteq V$. By \cite{Berkovich}, Proposition~2.5.20, $v$ belongs to the Shilov boundary for the strictly affinoid algebra $A_{U}$ of $U$. Therefore, the rational localization $A\to C\colonequals A_{U}$ defines a rational subset $\Spa(C, C^{+})$ of $\Spa(A, A^{\circ})$ as in (2).\end{proof}
We then have the following useful result, which was already observed by Poineau in \cite{Poineau13} and by Bhatt and Hansen in \cite{Bhatt-Hansen}.
\begin{prop}[Poineau \cite{Poineau13}, \S4, Bhatt-Hansen \cite{Bhatt-Hansen}, Proposition~2.9]\label{Abhyankar and Shilov}Let $A$ be an affinoid algebra in the sense of Tate over a nonarchimedean field $K$ and let $v$ be a continuous absolute value on $A$ over $K$. Then $v$ is a rational Abhyankar valuation if and only if it is weakly Shilov.\end{prop}
\begin{proof}If $v$ is weakly Shilov, Lemma~\ref{Weakly Shilov implies rationally weakly Shilov} implies that there exists a strictly affinoid neighborhood $V$ of $v$ in $\mathscr{M}(A)$ such that $v$ belongs to the Shilov boundary for the strictly affinoid algebra $A_{V}$. By \cite{Poineau13}, Lemme~4.4, applied to the strictly affinoid Berkovich-analytic space $V$, we see that $v$ is a rational Abhyankar valuation on $A_{V}$. But the completed residue field $\mathcal{H}(v)$ does not depend on whether we view $v$ as a point of $V=\mathscr{M}(A_{V})$ or as a point of $\mathscr{M}(A)$, so $v$ is a rational Abhyankar valuation on $A$. 

Conversely, suppose that $v$ is rational Abhyankar. Since $v$ is an absolute value on $A$, the underlying ring of $A$ is an integral domain, so the strictly affinoid space $\mathscr{M}(A)$ is irreducible. We can then apply \cite{Poineau13}, Corollaire~4.15, to conclude that $v$ is weakly Shilov. \end{proof}
The result \cite{Poineau13}, Corollaire~4.15, cited in the proof of Proposition~\ref{Abhyankar and Shilov} actually produces a strictly affinoid neighborhood $V$ of the rational Abhyankar valuation $v$ such that $v$ is the only point of the Shilov boundary for $A_{V}$. This type of argument already occurs in Berkovich's proof of \cite{Berkovich}, Proposition~2.4.4. Since this technique is going to be useful for us later when proving Theorem \ref{Topologically simple rules out Abhyankar}, we isolate the gist of it in the following lemma. 
\begin{lemma}\label{Localization of the Shilov boundary}Let $v$ be an element of the Shilov boundary for an affinoid algebra in the sense of Tate over a nonarchimedean field $K$. There exists an element $f\in A$ with $\vert f\vert_{\spc}=1$ such that $v\in \mathcal{M}(A\langle f^{-1}\rangle)$ and (the extension of) $v$ is the only element of the Shilov boundary for $A\langle f^{-1}\rangle$.\end{lemma}
\begin{proof}By Berkovich's description of the Shilov boundary for strictly affinoid algebras (\cite{Berkovich}, Proposition 2.4.4(iii)), the center $\spc(v)$ of $v$ on $A^{\circ}$ is a minimal open prime ideal of $A^{\circ}$, and we can again choose an element $\widetilde{f}\in\widetilde{A}$ which belongs to all minimal prime ideals of $\widetilde{A}$ other than $\mathfrak{p}=\spc(v)/A^{\circ\circ}$ (since $\widetilde{A}$ is Noetherian and, thus, has only finitely many minimal prime ideals, such an element can always be found). We then let $f\in A^{\circ}$ be a lift of $\widetilde{f}$, which necessarily satisfies $\vert f\vert_{\spc}=1$ since $\widetilde{f}\neq0$. By \cite{BGR}, Proposition 7.2.6/3, we have \begin{equation*}\widetilde{A\langle f^{-1}\rangle}=\widetilde{A}[\widetilde{f}^{-1}],\end{equation*}so $\mathfrak{p}\widetilde{A\langle f^{-1}\rangle}$ is the only minimal prime ideal of $\widetilde{A\langle f^{-1}\rangle}$. Since $A\langle f^{-1}\rangle$ is strictly affinoid, we can now apply Berkovich's \cite{Berkovich}, Proposition 2.4.4, to see that $v$ is the only element in the Shilov boundary for $A\langle f^{-1}\rangle$.\end{proof}
\begin{rmk}There is also a variant of the above lemma which applies to general affinoid algebras in the sense of Berkovich, but uses elements $f\in A$ with $\vert f\vert_{\spc}\neq 1$. Its proof uses Temkin's notion of graded reduction, see the proof of \cite{Temkin04}, Proposition 3.3.\end{rmk}
Another useful fact about Abhyankar valuations is the fact that the Abhyankar and rational Abhyankar conditions ascend and descend along finite extensions of affinoid algebras (Lemma~\ref{Finite extensions and Abhyankar valuations} below). To prove this, we use the following simple lemma on finite extensions of valued fields. 
\begin{lemma}\label{Finite field extensions and completion}Let $K\hookrightarrow L$ be a finite field extension, let $v$ be an absolute value on $L$ and $w=v\vert_{K}$. Then the completion $\widehat{L}^{v}$ of $L$ with respect to $v$ is also a finite extension of the completion $\widehat{K}^{w}$ of $K$ with respect to $w$.\end{lemma}
\begin{proof}Choose a basis $u_1,\dots, u_r$ of $L$ over $K$. Since $L$ is dense in $\widehat{L}^{v}$, the elements $u_1,\dots, u_r$ span a dense finite-dimensional $K$-linear subspace of $\widehat{L}^{v}$. But, by \cite{Schneider}, Proposition~4.13, there is only one norm topology on any finite-dimensional vector space over the nonarchimedean field $\widehat{K}^{w}$, and this norm topology is complete. It follows that the dense $\widehat{K}^{w}$-linear subspace of $\widehat{L}^{v}$ spanned by $u_1,\dots, u_r$ is already complete and thus coincides with $\widehat{L}^{v}$. In particular, $\widehat{L}^{v}$ is a finite field extension of $\widehat{K}^{w}$, as desired. \end{proof} 
\begin{lemma}\label{Finite extensions and Abhyankar valuations}Let $A\hookrightarrow B$ be a finite extension of affinoid algebras in the sense of Berkovich and let $v$ be a continuous absolute value on $B$ (over $K$). Then $v$ is Abhyankar (respectively, rational Abhyankar) if and only if the restriction $w$ of $v$ to $A$ is Abhyankar (respectively, rational Abhyankar).\end{lemma}
\begin{proof}The field extension $Q(A)\hookrightarrow Q(B)$ is finite, hence, by Lemma~\ref{Finite field extensions and completion}, so is the field extension $\mathcal{H}(w)\hookrightarrow \mathcal{H}(v)$. It follows that $\vert \mathcal{H}(v)^{\times}\vert/\vert\mathcal{H}(w)^{\times}\vert$ is torsion and $d_{K}(v)=d_{K}(w)$, and the assertion follows from Corollary~\ref{Abhyankar absolute values}. \end{proof}

\section{Topologically simple valuations and seminorms}\label{sec:topologically simple}

Recall that an element $x$ in a Hausdorff topological ring $A$ is called a topological unit (or topologically invertible) if there exists a net $(x_\lambda)_{\lambda}$ in $A$ such that $x_{\lambda}x \to 1$. This is equivalent to saying that the principal ideal generated by $x$ is a dense subset of $A$. Whenever $A$ has open unit group (for example, if $A$ is a Zariskian Huber ring), every topological unit is already a unit in $A$. However, this need not be true for more general topological rings. 
\begin{defn}[Topologically simple topological rings and modules]~\begin{enumerate}[(1)]\item A topological module $M$ over a topological ring $A$ is called topologically simple if the closure of the zero submodule is the only closed submodule of $M$.\item A topological ring $A$ is said to be topologically simple if it is topologically simple as a topological module over itself, i.e. if its only closed ideal is the closure of the zero ideal.\end{enumerate}\end{defn}
\begin{rmk} The term 'topologically simple'  is meant to bring to mind an obvious analogy with the notion of a simple module in algebra. It is also commonly used in the theory of Banach algebras without unit.
\end{rmk}

\begin{lemma}
    \label{Topologically simple topologies and topological units}A topological ring $A$ is
topologically simple if and only if every nonzero element $f$ of the quotient $\overline{A}=A/\overline{(0)}$ by the closure of the zero ideal is a topological unit.\end{lemma}
\begin{proof}It is clear that $A$ is topologically simple if and only if $\overline{A}$ is topologically simple. Hence it suffices to assume that $A$ is Hausdorff. Since the closure of an ideal in a topological ring is always a (not necessarily proper) ideal, $A$ being topological simple is equivalent to every ideal of $A$ being dense. This, in turn, is equivalent to every principal ideal being dense and, since $A$ is Hausdorff, this is the same every element of $A$ being a topological unit. \end{proof}
We can also talk about valuations and seminorms being topologically simple. Since we want to include both the case of a valuation which is not necessarily of rank one and of a seminorm which is not necessarily multiplicative, we phrase our definition in terms of generalized seminorms, in the sense of Paugam \cite{Paugam09}.
\begin{defn}[Pseudo-uniformizer]\label{Pseudo-uniformizer}Let $v: A\to \Gamma$ be a nonarchimedean generalized seminorm in the sense of \cite[Definition~14]{Paugam09} (in particular, $v$ could be a valuation of arbitrary rank, or an arbitrary seminorm $A\to \mathbb{R}_{\geq0}$). An element $\varpi\in A$ is called a pseudo-uniformizer for $v$ is a unit in $A$, if $\varpi$ is multiplicative with respect to $v$ (i.e., $v(\varpi f)=v(\varpi)v(f)$ for all $f\in A$) and if $0<v(\varpi)<1$.\end{defn}
Note that this implies that $v(\varpi^{-1})=v(\varpi)^{-1}>1$, so $\varpi$ is not a unit in the closed unit ball $A_{v\leq1}$ of $v$.    
\begin{defn}[Topologically simple valuations and seminorms]
\label{Topologically simple valuation}Let $v: A\to \Gamma$ be a nonarchimedean generalized seminorm in the sense of \cite[Definition~14]{Paugam09}. Then $v$ is said to be topologically simple if there exists a pseudo-uniformizer $\varpi$ of $v$ such that the Tate ring $A$ with pair of definition $(A_{v\leq1}, \varpi)$ is a topologically simple topological ring.\end{defn}
\begin{rmk}One could imagine that, for some higher rank valuation $v$, the condition in Definition \ref{Topologically simple valuation} could be satisfied for some pseudo-uniformizer, but not every pseudo-uniformizer $\varpi$ of $v$. We will show in Theorem \ref{Pruefer-Manis vs. topologically simple} that this cannot happen.\end{rmk} 
It is readily seen that a Hausdorff topological ring $A$ is topologically simple if and only if
every ideal (except the zero ideal) is dense in $A$, and this occurs if and only if every non-zero
element of $A$ is a topological unit. Such a ring $A$ clearly must be an integral domain and one can
ask whether any topologically simple $A$ must already be a field. The answer is known to be negative
if $A$ is a topological algebra over $\mathbb{C}$. Indeed, Atzmon \cite{Atzmon} has constructed a
complete metrizable topological algebra $A$ over $\mathbb{C}$ such that every non-zero element of
$A$ is a topological unit but $A$ is not a field. However, we note that, the algebra $A$ from
\cite{Atzmon} being complete, its topology cannot be defined by a norm. Thus we may modify our
question to ask whether there exists a topologically simple normed ring $A$ which is not a field. By
the lemma stated below it suffices to ask for a normed ring $A$ such that $A$ has no closed prime
ideals or such that every element of the Gelfand spectrum of $A$ is an absolute value. 
\begin{lemma}[cf.~\cite{Dine22}, Lemma 3.12] \label{Every closed ideal is contained in a closed
    prime ideal} Let $A$ be a seminormed ring. Then every ideal $I \subsetneq A$ which is not dense
    in $A$ is contained in the kernel of some element of the Gelfand spectrum of $A$. In particular, $I$ is contained in a spectrally reduced prime ideal
\end{lemma}
\begin{proof}Follows from the non-emptiness of the Gelfand spectrum of the quotient seminormed ring $A/I$. \end{proof} 
Every element $\phi$ in the Gelfand spectrum of a topologically simple normed ring is an absolute value (otherwise $\ker(\phi)$ would be a closed ideal). Hence \cite[Theorem~1.4]{Escassut1}, implies that a topologically simple uniform normed algebra over a nonarchimedean field $K$ has no non-zero topological divisors of zero. We continue with some further simple observations concerning the properties of topologically simple rings. 

\begin{lemma}\label{Torsion modules}Let $A$ be a topologically simple Hausdorff topological ring and let $M$ be an $A$-module. Then either $M$ is a torsion-free module or $M$ admits no Hausdorff $A$-module topology. More precisely, any $A$-module topology on a torsion $A$-module coincides with the indiscrete topology. \end{lemma} 
\begin{proof}If $M$ is a topological $A$-module and the subspace topology on every finitely generated submodule of $M$ coincides with the indiscrete topology, then the topology on $M$ is the indiscrete topology (for a non-empty open subset $U$ of $M$, pick any $x\in U$, then $U\cap (x)_{A}$ is an open subset of the finitely generated submodule $(x)_{A}$ of $M$, so, by assumption, $0\in U\cap (x)_{A}\subseteq U$). Hence it suffices to prove that any $A$-module topology on a finitely generated torsion $A$-module coincides with the indiscrete topology. Let $M$ be a finitely generated torsion $A$-module. Then the annihilator~$\Ann(M)$ of $M$ is a non-zero ideal of $A$. Endow $M$ with any $A$-module topology. Then the map $A \times M \to M$ becomes continuous (when $A \times M$ is given the product topology). Since $\Ann(M) \subseteq{A}$ is dense, $\Ann(M) \times M$ is dense in $A \times M$ and, by continuity of $A \times M \to M$, we see that $(0)=\Ann(M)M$ is dense in $M$. That is, the topology on $M$ is indiscrete.\end{proof}
Recall that a module over a ring $A$ is said to be simple if it contains no non-zero proper submodules. Any simple $A$-module is isomorphic to $A/\mathfrak{m}$ for $\mathfrak{m} \subsetneq A$ a maximal ideal. For $M$ an $A$-module, recall that a composition series for $M$ a finite chain of $A$-submodules
\begin{equation*}0=M_{0} \subseteq M_{1} \subseteq \dots \subseteq M_{n}=M
\end{equation*}such that for all $i\geq 1$ the quotient $M_{i}/M_{i-1}$ is a simple $A$-module. By the Jordan-Hölder Theorem, a module $M$ has a composition series if and only if it is of finite length.   

\begin{lemma} Let $A$ be a topological ring, not a field, such that every maximal ideal of $A$ is dense (in particular, we could take $A$ topologically simple Hausdorff). For every $A$-module $M$ of finite length, any $A$-module topology on $M$ coincides with the indiscrete topology. In particular, if $A$ is an Artinian ring which is not a field, then the only topology turning $A$ into a topologically simple topological ring is the indiscrete topology.\end{lemma}

\begin{proof} 
    Consider a composition series 
    \begin{equation*}
        0=M_{0} \subseteq M_{1} \subseteq \dots \subseteq M_{n} = M
    \end{equation*} 
    for $M$. For every $i \geq 1$ the quotient $M_{i}/M_{i-1}$ is isomorphic to $A/ \mathfrak{m}$ for some maximal ideal $\mathfrak{m} \subsetneq{A}$. Since $\mathfrak{m}$ is dense in $A$, the natural topology on the quotient $A/\mathfrak{m}$ coincides with the indiscrete topology. But the natural topology on an $A$-module is the finest $A$-module topology, so any $A$-module topology on $M_{i}/M_{i-1}$ coincides with the indiscrete topology. Now endow $M$ with some $A$-module topology. By the above, $M_{i-1}$ is dense in $M_{i}$ for all $i \geq 1$. Then $(0)$ is dense in $M$.
\end{proof}

The next lemma involves the notion of Zariskisation for Huber rings discussed in the \cite{Tanaka} and \cite{Dine22}, \S3. 
\begin{lemma}
    \label{Topologically simple rings and Zariskisation}
    Let $A$ be a Huber ring. 
    Then the following are equivalent:\begin{enumerate}[(1)]\item $A$ is topologically simple.\item The Zariskisation $\overline{A}^{\Zar}$ of the Hausdorff quotient $\overline{A}=A/\overline{(0)}$ is a field.\end{enumerate}\end{lemma}
\begin{proof}It suffices to assume that $A$ is Hausdorff. If $A^{\Zar}$ is a field, then every element of $A$ is topologically invertible since $A$ is a dense subring of $A^{\Zar}$. Now suppose that $A$ is a Hausdorff Huber ring with every non-zero element of $A$ topologically invertible. Since $A^{\Zar}$ is a localization of $A$, any fixed non-zero prime ideal $\mathfrak{p}$ of $A^{\Zar}$ pulls back to a non-zero prime ideal in $A$. In particular, $\mathfrak{p}$ contains a non-zero element $f$ that is an image of a non-zero element of $A$. Since the map $A \to A^{\Zar}$ is continuous, $f$ is a topological unit in $A^{\Zar}$. But, by definition of a Zariskian Huber ring, the group of units of $A^{\Zar}$ is open in $A^{\Zar}$, so $f$ is in fact a unit in $A^{\Zar}$ contradicting the assumption that $\mathfrak{p}$ is a (proper) prime ideal. Thus there exist no non-zero prime ideals in $A^{\Zar}$, i.e. $A^{\Zar}$ is a field.
\end{proof}
\begin{lem}
    \label{lem:topsimp:complete}
    Let \(A\) be a topological ring. If the completion of $A$ is a field, then $A$ is topologically simple. If $A$ is a Huber ring and its topology is defined by an absolute value, then the converse is also true.
\end{lem}
\begin{proof}
    Suppose the completion of \(A\) is a field. 
    Pick \(f \in A\), then there exists an element \(g \in \widehat{A}\) such that
    \(fg = 1\) in \(\widehat{A}\). 
    Then choose a Cauchy sequence \(g_{n} \xrightarrow{} g\) 
    in \(A\). 
    Then \(g_{n}f \xrightarrow{} 1\), so \(f\) is a topological unit.
    Conclude by Lemma \ref{Topologically simple topologies and topological units} that $A$ is topologically simple.
    
    Conversely, let \(A\) be a topologically simple Huber ring whose topology is defined by an absolute value.
    Then the completion of \(A\) factors through the Zariskization \(A^{\Zar}\), which is a field by Lemma \ref{Topologically simple rings and Zariskisation}.
    Since \(A\) is Huber and its topology is defined by an absolute value, the completion \(\widehat{A}\) of $A$ 
    is also the completion of \(A^{\Zar}\) with respect to an absolute value,
    and the completion of any field with respect to an absolute value is a field.
\end{proof}
\begin{lemma}\label{Topologically simple norms and rational subsets}Let $A$ be a Huber ring and let $v$ be a topologically simple continuous rank one valuation on $A$. If $A\to B$ is a continuous map of Huber rings such that $B$ has a dense $A$-subalgebra which is a localization of $A$ and such that $v$ extends to a continuous valuation $w$ on $B$, then $w$ is also topologically simple. \end{lemma}
\begin{proof}Denote by $\widehat{A}^{v}$, $\widehat{B}^{w}$ the completion of $A$, $B$ with respect to $v$, $w$. By continuity, the map $A\to B$ extends to a continuous map $\widehat{A}^{v}\to \widehat{B}^{w}$. By assumption, there exists a multiplicative subset $S$ of $A$ such that $S^{-1}A$ is a dense subring of $B$. Then the image of $S^{-1}A$ in $\widehat{B}^{w}$ is also a dense subring of $\widehat{B}^{w}$ (the valuation $w$ on $B$ being continuous). A fortiori, the image of $S^{-1}\widehat{A}^{v}$ in $\widehat{B}^{w}$ is dense in $\widehat{B}^{w}$. But, by virtue of Lemma \ref{lem:topsimp:complete}, $\widehat{A}^{v}$ is a field, so $S^{-1}\widehat{A}^{v}=\widehat{A}^{v}$. It follows that the map $\widehat{A}^{v}\to \widehat{B}^{w}$ has dense image. Since $w$ restricts to $v$ on $A$, we conclude that $\widehat{A}^{v}=\widehat{B}^{w}$. In particular, $\widehat{B}^{w}$ is a field and thus $w$ is topologically simple. \end{proof}
\begin{rmk}In particular, the above lemma applies to rational localizations $A\to B$ of a Huber ring $A$.\end{rmk}

In the following we will construct an example (or rather a class of examples) of a topologically simple normed ring which is far from being a field algebraically. 

\begin{thm} \label{Topologically simple norms on the polynomial algebra} For $K$ a nonarchimedean field the following statements are equivalent:\begin{enumerate}[(1)]\item The completed algebraic closure $\widehat{\overline{K}}$ is spherically complete.\item There exists no topologically simple $K$-algebra norm $\vert \cdot\vert_{A}$ on the $K$-algebra $A=K[T]$.\end{enumerate}More precisely, for any nonarchimedean field $K$, the topologically simple absolute values on $K[T]$ extending the given absolute value on $K$ are precisely the type (IV) points. \end{thm}
\begin{proof}Assume that $\widehat{\overline{K}}$ is not spherically complete. By Lemma \ref{Every closed ideal is contained in a closed prime ideal}, every closed ideal in a normed ring is contained in a closed prime ideal, so ($A$ being of Krull dimension $1$) it suffices to find $\vert\cdot\vert_{A}$ such that no maximal ideal of $A$ is closed in the topology induced by this norm. Since $\widehat{\overline{K}}$ is not spherically complete, we can find a nested sequence of closed disks $B(a_1,r_1) \supsetneq B(a_2,r_2) \supsetneq B(a_3,r_3) \supsetneq \dots $ in $\widehat{\overline{K}}$ having empty intersection. Consider the type $(IV)$ point corresponding to this nested sequence of disks, i.e., the multiplicative seminorm\begin{equation*} \vert f\vert_{A}:=\inf_{n} \vert f\vert_{B(a_n,r_n)} \ \ (f \in A).
\end{equation*}Let $\mathfrak{m} \subsetneq{A}$ be a maximal ideal. By Hilbert's Nullstellensatz, $\mathfrak{m}$ is of the form $\ker(\vert \cdot\vert_{a})$ for some point $a \in \overline{K}$ with $\vert f\vert_{a}=\vert f(a)\vert$ ($f \in A$). The maximal ideal $\mathfrak{m}$ is closed if and only if the seminorm $\vert\cdot\vert_{a}$ is bounded by $\vert \cdot\vert_{A}$. But in this case $\vert\cdot\vert_{a} \leq \vert\cdot\vert_{B(a_n,r_n)}$ for every $n$. In particular, for every $n$ we have $\vert a-a_n\vert=\vert T-a_n\vert_{a} \leq \sup_{b \in B(a_n,r_n)}\vert T-a_n\vert_{b}=\sup_{b \in B(a_n,r_n)}\vert b - a_n\vert \leq r_n$. That is, $a \in \bigcap_{n} B(a_n,r_n)$, a contradiction.

Conversely, assume that the completed algebraic closure $\widehat{\overline{K}}$ of the ground field
$K$ is spherically complete and let $\vert\cdot\vert_{A}$ be a norm on $A=K[T]$ restricting to the
given absolute value on $K$. We want to show that $\vert\cdot\vert_{A}$ is not topologically simple.
If $\vert\cdot\vert_{A}$ is topologically simple, then every element in the Gelfand spectrum of the normed ring $(A, \vert\cdot\vert_{A})$ is an absolute value. Besides, if $\vert\cdot\vert_{A}$ is topologically simple, then every absolute value $\vert\cdot\vert$ bounded above by $\vert\cdot\vert_{A}$ is topologically simple as well (an ideal which is closed with respect to $\vert\cdot\vert$ is also closed with respect to the stronger norm $\vert\cdot\vert_{A}$). Hence it suffices to prove our claim in the case when $\vert\cdot\vert_{A}$ is an absolute value, i.e., when the norm $\vert\cdot\vert_{A}$ defines a point on the Berkovich affine line $\mathbb{A}_{K}^{1,\an}$. Since $\widehat{\overline{K}}$ is spherically complete, $\lVert\cdot \rVert$ must then be of the form\begin{equation*}\vert f\vert_{A}=\vert f\vert_{B(a,r)} \ \ (f \in A)\end{equation*}for some closed disk $B(a,r)$ in $\widehat{\overline{K}}$. Consider the multiplicative seminorm $\vert\cdot\vert_{a}$ on $K[T]$ given by $\vert f\vert_{a}=\vert f(a)\vert$ for $f \in A$ (evaluation at the centre of the disk). Then clearly $\vert \cdot\vert_{a} \leq \vert\cdot\vert_{A}$, so $\vert\cdot\vert_{a}$ is continuous with respect to the norm $\vert\cdot\vert_{A}$. Hence the maximal ideal $\ker(\vert\cdot\vert_{a})$ is closed with respect to the the topology defined by $\vert\cdot\vert_{A}$. In particular, the norm $\vert\cdot\vert_{A}$ is not topologically simple.\end{proof}
We also have the following variant of the theorem for any closed disk. 

\begin{thm}\label{Topologically simple norms on the Tate algebra}For $K$ a nonarchimedean field, the following statements are equivalent:\begin{enumerate}[(1)]\item The completed algebraic closure $\widehat{\overline{K}}$ of $K$ is spherically complete.\item There exists no topologically simple norm $\vert \cdot\vert$ on $K\langle T\rangle$ bounded above by the Gauss norm. \item There exists $r>0$ such that there is no topologically simple norm $\vert \cdot\vert$ on $K\langle r^{-1}T\rangle$ bounded above by the Gauss norm. \item For all $r>0$, there exists no topologically simple norm $\vert\cdot\vert$ on $K\langle r^{-1}T\rangle$ bounded above by the Gauss norm.\end{enumerate}More precisely, for any nonarchimedean field $K$ and any radius $r>0$, the topologically simple absolute values on $K\langle r^{-1}T\rangle$ bounded above by the Gauss norm are precisely the type $(IV)$ points in the closed ball $B(0, r)$. \end{thm}
\begin{proof}The equivalence of the first and last conditions is exactly Theorem \ref{Topologically simple norms on the polynomial algebra}, since the Berkovich affine line $\mathbb{A}_{K}^{1,\an}$ can be written as a union of Berkovich disks $\mathscr{M}(K\langle r^{-1}T\rangle)$ of varying radius $r>0$. In particular, condition (1) implies all of the other conditions. Thus it remains to show that, if $\widehat{\overline{K}}$ is not spherically complete and $r>0$ is any radius, there exists a topologically simple norm $\vert\cdot\vert$ on $K\langle r^{-1}T\rangle$ bounded above by the Gauss norm. But for any nested sequence of disks in $\widehat{\overline{K}}$ with empty intersection, we can shift centers of the disks, and possibly pass to a subsequence, to obtain an analogous nested sequence of disks which are contained in the disk $B(0, r)$. But then the corresponding type $(IV)$ point $\vert\cdot\vert$ becomes bounded by the Gauss norm on $K\langle r^{-1}T\rangle$. This concludes the proof since we have already seen in Theorem \ref{Topologically simple norms on the polynomial algebra} that every type $(IV)$ point is a topologically simple absolute value.\end{proof}

Using Theorem \ref{Topologically simple norms on the polynomial algebra}, we can now describe all
completions of the \(n\)-variable polynomial ring (or the \(n\)-variable Tate algebra) with respect to topologically simple absolute values.

\begin{prop}
    \label{prop:topsimp:nvar}
    Let \(K\) be a nonarchimedean field, and let \(n \in N\),
    and let \(v\) be a topologically simple absolute value on 
    the polynomial ring \(A = K[T_{1}, \ldots, T_{n}]\).
    Let \(w\) denote the restriction of \(v\) to 
    \(K[T_{1}, \ldots, T_{n-1}]\).
    Then, the canonical extension o \(v\) to
    \(\mathcal{H}(w)[T_{n}]\) is 
    a type \((IV)\) point.
    Moreover, this holds for any permutation of the variables
    \(T_{1}, \ldots, T_{n}\).

    Similarly, for \(A = K\langle T_{1}, \ldots, T_{n}\rangle \),
    \(v\) topologically simple on \(A\),
    and \(w\) the restriction
    to \(K\langle T_{1}, \ldots, T_{n-1}\rangle \),
    we have that the canonical extension of \(v\) 
    to \(\mathcal{H}(w)\langle T_{n}\rangle \) 
    is a type (IV) point.
    Moreover, this holds for any permutation of the variables
    \(T_{1}, \ldots, T_{n}\).
\end{prop}

\begin{proof}
    In the first case, \(A\) is dense in its completion, and thus so is
    \(\mathcal{H}(w)[T_{n}]\) (by the universal 
    properties of localization and completion 
    and Lemma \ref{lem:topsimp:complete}).
    Then, by Theorem \ref{Topologically simple norms on the polynomial algebra},
    \(v\) must be a type (IV) point.
    We may permute the variables first and the argument remains unchanged.
    To prove the statement for the Tate algebra,
    replace Theorem \ref{Topologically simple norms on the polynomial algebra}
    with Theorem \ref{Topologically simple norms on the Tate algebra}.
\end{proof}

\subsection{Topologically simple valuations and Pr\"ufer-Manis valuations}

We conclude this section by relating the notion of topologically simple valuations to the notion of so-called Prüfer-Manis valuations studied in the book \cite{Knebusch-Zhang} of Knebusch and Zhang. This relationship will also reveal that the property of a (possibly higher-rank) valuation being topologically simple does not depend on the choice of a pseudo-uniformizer $\varpi$.
\begin{defn}
    Let \(v: R\to \Gamma\cup\{0\}\) be a valuation on \(R\). 
    Then \(v\) is \textit{Manis} if 
    \(v(R^{\times}) = c\Gamma_{v}\), i.e., if the value semigroup $v(R\setminus\{0\})$ of $v$ is a subgroup of $\Gamma$. 
\end{defn}
\begin{defn}
    Let \(R\) be a ring, 
    and let \(A\) be a subring of \(R\) and
    \(\mathfrak{p}\) be a prime ideal of \(A\).
    Then \((A,\mathfrak{p})\) is a 
    \textit{Manis pair} if
    \((A, \mathfrak{p})=(A_{v\leq1}, A_{v<1})\)
    for some Manis valuation \(v\).\end{defn}
\begin{defn}
    Let \(A \ins R\) be an extension of rings.
    The extension is called a Pr\"ufer
    extension, or $A$ is called a Prüfer subring of $R$, if 
    \((A_{\mathfrak{m}} \mathfrak{m}A_{\mathfrak{m}})\) is a Manis pair in 
    \(R_{\mathfrak{m}}\) for
    all maximal ideals \(\mathfrak{m}\) of \(A\).\end{defn} 
\begin{defn}[Prüfer-Manis valuation]A valuation $v: A\to\Gamma\cup\{0\}$ on a ring $R$ is called Prüfer-Manis if $v$ is a Manis valuation and the closed unit ball $R_{v\leq1}$ is a Prüfer subring of $R$. A Prüfer-Manis subring of $R$ is a subring of the form $A_{v\leq1}$ for some Prüfer-Manis valuation $v$ on $R$, i.e., it is a Prüfer subring which is the closed unit ball of a Manis valuation $v$ on $R$. \end{defn}
We recall the following result of Knebusch-Zhang.
\begin{prop}[\cite{Knebusch-Zhang}, Ch.~I, Prop.~5.1(iii)]\label{Pruefer vs. Pruefer-Manis}Let $v$ be a valuation on a ring $R$. Then $R_{v\leq1}$ is a Prüfer-Manis subring of $R$ if and only if it is a Prüfer-Manis subring of $R$. In other words, if a Prüfer subring $R_{0}$ of $R$ is also the closed unit ball of a valuation $v: R\to\Gamma\cup\{0\}$ (in the terminology of \cite{Knebusch-Zhang}, $R_{0}$ is a valuation subring of $R$), then the valuation $v$ must automatically be Manis.\end{prop}  
We also recall the following observation from \cite{Dine25} which relates the notion of Prüfer subrings studied by Knebusch-Zhang to the notion of a $\varpi$-valuative ring considered in the book \cite{FK} of Fujiwara-Kato. 
\begin{prop}[\cite{Dine22}, Proposition 4.20]\label{Pruefer subrings and valuative rings}Let $\varpi$ be a non-zero-divisor in a ring $A$. Then $A$ is a Prüfer subring of $A[\varpi^{-1}]$ if and only if $A$ is a $\varpi$-valuative ring in the sense of Fujiwara-Kato \cite{FK}, i.e., every finitely generated ideal $I$ of $A$ which contains a power of $\varpi$ is invertible.\end{prop}
One of the useful properties of Prüfer subrings proved in the book of Knebusch-Zhang is their stability under taking quotients.
\begin{prop}[\cite{Knebusch-Zhang}, Ch.~I, Prop.~5.7-5.8]\label{Pruefer subrings and quotients}Let $A$ be a Prüfer subring of a ring $R$ and let $\varphi: R\to S$ be a ring map. Then $\varphi(A)$ is a Prüfer subring of $\varphi(R)$. If $I$ is an ideal of a ring $R$ contained in a subring $A$, then $A$ is a Prüfer subring of $R$ if and only if $A/I$ is a Prüfer subring of $R/I$.\end{prop}
\begin{cor}\label{Pruefer-Manis valuations and Hausdorff quotient}Let $v$ be a valuation on a ring $R$. Then $v$ is Prüfer-Manis if and only if the induced valuation $\overline{v}$ on the quotient $\overline{R}=R/\ker(v)$ is Prüfer-Manis.\end{cor}
\begin{proof}By definition of the induced valuation on $\overline{R}$, we have \begin{equation*}\overline{R}_{\overline{v}\leq1}=R_{v\leq1}/\ker(v).\end{equation*}Since $\ker(v)$ is an ideal of $R$ contained in $R_{v\leq1}$, Proposition \ref{Pruefer subrings and quotients} tells us that $R_{v\leq1}$ is a Prüfer subring of $R$ if and only if $\overline{R}_{\overline{v}\leq1}$ is a Prüfer subring of $\overline{R}$. We conclude by applying Proposition \ref{Pruefer vs. Pruefer-Manis}.\end{proof}
Finally, we recall the following theorem of Knebusch-Zhang.
\begin{thm}[\cite{Knebusch-Zhang}, Ch.~III, Theorem 1.8]\label{Maximal ideals in a Pruefer-Manis subring}Let $A$ be a Prüfer-Manis subring of a ring $R$. Then $A$ has only one maximal ideal $\mathfrak{m}$ such that $\mathfrak{m}R=R$. If $v$ is a Prüfer-Manis valuation such that $A=R_{v\leq1}$, then $\mathfrak{m}=R_{v<1}$.\end{thm} 
\begin{thm}\label{Pruefer-Manis vs. topologically simple}Let $v: R\to\Gamma\cup\{0\}$ be a non-trivial valuation on a ring $R$. Then $v$ is Prüfer-Manis if and only if it is topologically simple. Moreover, in this case, the Tate ring $R$ with pair of definition $(R_{\leq1}, \varpi)$ is topologically simple for every choice of pseudo-uniformizer $\varpi$ of $v$.  \end{thm}
\begin{proof}Set $A=R_{v\leq1}$. Let $\varpi\in A$ be a pseudo-uniformizer of $v$. By Corollary \ref{Pruefer-Manis valuations and Hausdorff quotient}, we may assume that $\ker(v)=0$.

Suppose that $v$ is Prüfer-Manis. Note that $R=A[\varpi^{-1}]$. By Theorem \ref{Maximal ideals in a
Pruefer-Manis subring}, $R_{v<1}$ is the unique maximal ideal of $A$ containing $\varpi$. Hence, if
$\mathfrak{p}$ is a prime ideal of $A$, then either $\varpi\in\mathfrak{p}$ and $\mathfrak{p}$ is
contained in $R_{v<1}$, or $(\mathfrak{p}, \varpi)_{A}=A$, in which case $\mathfrak{p}\cap
(1+(\varpi)_{A})\neq\emptyset$. It follows that the $\varpi$-adic Zariskisation \begin{equation*}A^{\Zar}=(1+(\varpi)_{A})^{-1}A\end{equation*}is a local ring. By Proposition \ref{Pruefer subrings and valuative rings}, $A$ is a $\varpi$-valuative ring, so $A^{\Zar}$ is also a $\varpi$-valuative ring, being a localization of $A$ (\cite{FK}, Ch.~0, Proposition 8.7.2). Since $\ker(v)=0$, we know that $A$ is $\varpi$-adically separated. Then $A^{\Zar}$ is also $\varpi$-adically separated (\cite{Tanaka}, Lemma 3.19(2)). By \cite{Dine25}, Proposition 4.23, it then follows that $A^{\Zar}$ is $\varpi$-adically separated valuation ring and thus $A^{\Zar}[\varpi^{-1}]$ is a field. By \cite{Tanaka}, Lemma 3.5, $A^{\Zar}[\varpi^{-1}]$ is precisely the Zariskisation $R^{\Zar}$ of the Tate ring $R$, with pair of definition $(A, \varpi)$. We conclude that $R^{\Zar}$ is a field and thus that the valuation $v$ is topologically simple, by Lemma \ref{lem:topsimp:complete}

Conversely, suppose that $v$ is topologically simple. Then, by Lemma \ref{lem:topsimp:complete}, $R^{\Zar}$ is a field and thus $A^{\Zar}$ (which coincides with the closed unit ball of the canonical extension of $v$ to $R^{\Zar}$) is a valuation ring. In particular, $A^{\Zar}$ is local. Since $A^{\Zar}/\varpi=A/\varpi$, this implies that $A$ has only one maximal ideal $\mathfrak{m}$ containing $\varpi$. Since $1+(\varpi)_{A}\subseteq A\setminus\mathfrak{m}$, the map $A\to A_{\mathfrak{m}}$ factors through $A\to A^{\Zar}$; therefore, $A_{\mathfrak{m}}$ is a localization of the valuation ring $A^{\Zar}$ and thus itself a valuation ring. In particular, $A_{\mathfrak{m}}$ is a $\varpi$-valuative ring. But for any other maximal ideal $\mathfrak{m}'\neq\mathfrak{m}$ of $A$, we know that $\varpi$ is a unit in $A_{\mathfrak{m}'}$, so $A_{\mathfrak{m}'}$ is (trivially) also a $\varpi$-valuative ring. By \cite{FK}, Ch.~0, Proposition 8.7.3, we see that $A$ is itself $\varpi$-valuative. By Proposition \ref{Pruefer subrings and valuative rings}, this means that $A$ is a Prüfer subring of $R$, in which case the valuation $v$ with $A=R_{v\leq1}$ is Prüfer-Manis, by Proposition \ref{Pruefer vs. Pruefer-Manis}.\end{proof}

\subsection{Relation to Abhyankar valuations}

We now prove Theorem \ref{thm:ab:topsimp} from the introduction, i.e., we show that two of the conditions on a continuous absolute value on an affinoid
algebra $A$ we have studied so far, namely, ``rational Abhyankar'' and ``topologically simple'', are mutually
exclusive for affinoid algebras in the sense of Tate (i.e., strictly affinoid algebras). 
Before we start, let's exposit an example that makes such a statement rather subtle.



\begin{exam}In the same vein, consider now the strictly affinoid algebra 
    \begin{equation*}A = K\langle X\rangle \langle s / X\rangle \colonequals K\langle X, T / s\rangle / (X
    (T / s) - 1),\end{equation*}where now 
    \(s \in \sqrt{\vert K^{\times}\vert}\) and \(s < 1\), that is, the affinoid algebra of the annulus with exterior radius $1$ and rational interior radius. Consider some number $r\in (s, 1)$ such that $r\not\in \sqrt{\vert K^{\times}\vert}$.
    Then the $r$-Gauss norm $v$ on $A$ dominates only those points of $\mathscr{M}(A)$ which lie
    on the irrational circle $\{\, \vert X\vert=r\,\}$:
    Indeed, the points $w\in \mathscr{M}(A)$ with $w(X)>r$ satisfy $w(X)>v(X)$ while the points $w$ with
    $w(X)<r$ satisfy $w(X^{-1})>v(X^{-1})$.
    But the circle $\{\, \vert X\vert=r\,\}$ consists of only one point, namely, the $r$-Gauss norm itself.
    It follows that
    the Gelfand spectrum of the normed ring $(A, v)$ consists of a single point, so the completion
    $\widehat{A}^{v}$ of $A$ with respect to $v$ must be a nonarchimedean field
    by \cite{Kedlaya18}, Lemma 2.4.
    (indeed, it is once again equal to $K_r$). By Lemma \ref{lem:topsimp:complete}, $v$
    is a topologically simple valuation on $A$. But $v$ is also Abhyankar, since $\mathcal{H}(v)=K_{r}$
    satisfies $d_{K}(K_{r})=1=\dim(A)$.
    \end{exam}
    Our next result, announced as Theorem \ref{thm:ab:topsimp} in the introduction, states that such phenomena as explained in the above two examples cannot occur for \textit{rational} Abhyankar valuations on strictly affinoid algebras.
\begin{thm}\label{Topologically simple rules out Abhyankar}Let $K$ be a nonarchimedean field, let $A$ be a $K$-affinoid algebra in the sense of Berkovich and let $v$ be a continuous absolute value on $A$ over $K$ which belongs to some strictly affinoid domain of $\mathscr{M}(A)$. If $v$ is topologically simple and $A$ is not a finite field extension of $K$, then $v$ is not a rational Abhyankar valuation.\end{thm}
We prove the theorem by assuming that $v$ is rational Abhyankar and using Theorem \ref{Abhyankar and Shilov} as well as the technique of ``Shilov localization", Lemma \ref{Localization of the Shilov boundary}, to pass to a strictly affinoid domain $V\subseteq \mathscr{M}(A)$ such that $v$ is equal to the spectral norm on $A_{V}$. We then leverage the topologically simple condition to deduce that $A$ must be a finite field extension of $K$.

To carry out this argument, we need to know that the topologically simple condition is stable under passage to affinoid subdomains. Note that it already follows from Lemma \ref{Topologically simple norms and rational subsets} that a topologically simple continuous absolute value $v$ on $A$ remains topologically simple after passage to a rational subdomain of $\mathscr{M}(A)$ containing $v$. Indeed, in this case, the affinoid algebra $A_{V}$ is a completed algebraic localization of $A$, by \cite{Berkovich}, Corollary 2.2.10, and then Lemma \ref{Topologically simple norms and rational subsets} applies. We now prove a lemma which extends this observation to more general affinoid domains of $\mathscr{M}(A)$.
\begin{lemma}\label{Topologically simple absolute values and affinoid domains}Let $A$ be a $K$-affinoid algebra in the sense of Berkovich over some nonarchimedean field $K$, let $v$ be a continuous absolute value on $A$ over $K$ and let $V$ be an affinoid domain with $v\in V$. If $v$ is topologically simple, then so is the extension of $v$ to $A_{V}$.\end{lemma}
\begin{proof}Since $v$ is an absolute value, $A$ is an integral domain. It then follows from \cite{Berkovich}, Proposition 2.3.4, that the afffinoid algebra $A_{V}$ is reduced. By \cite{Berkovich}, Proposition 2.1.4, this means that the spectral seminorm on $A_{V}$ is a norm, so $v$ is an absolute on $A_{V}$ (i.e., $\ker(v)=0$). By Lemma \ref{lem:topsimp:complete}, the topologically simple condition means that the completion $\widehat{A}^{v}$ of $A$ with respect to $v$ is a field. In particular, $\widehat{A}^{v}$ must contain the fraction field $Q(A)$ of $A$. Since already $A$ is dense in $\widehat{A}^{v}$, so is $Q(A)$. It follows that \begin{equation*}\widehat{A}^{v}=\mathcal{H}(v),\end{equation*}the completed residue field of $\mathscr{M}(A)$ at $v$. In other words, $A$ is dense in $\mathcal{H}(v)$. By the universal property of affinoid domains, the map $A\hookrightarrow \mathcal{H}(v)$ factors through $B\to \mathcal{H}(v)$, where the latter map is also injective since we have seen that $v$ has zero kernel on $A_{V}$. It follows that $A_{V}$ is also dense in $\mathcal{H}(v)$, that is, the completion of $A_{V}$ with respect to $v$ is again the nonarchimedean field $\mathcal{H}(v)$. We conclude by applying Lemma \ref{lem:topsimp:complete}. \end{proof}  
\begin{proof}[Proof of Theorem \ref{Topologically simple rules out Abhyankar}]Suppose that $v$ is both rational Abhyankar and topologically simple. By Proposition \ref{Abhyankar and Shilov}, there exists a strictly affinoid domain $V$ of $\mathscr{M}(A)$ containing $v$ such that $v$ belongs to the Shilov boundary for $A_{V}$. By Lemma \ref{Localization of the Shilov boundary}, up to shrinking $V$, we may assume that $v$ is actually the only point in the Shilov boundary for $A_{V}$, i.e., it is equal to the spectral seminorm on $A_{V}$. By Lemma \ref{Topologically simple absolute values and affinoid domains}, this means that $v$ is also topologically simple as a valuation on $A_{V}$. On the other hand, $A$ being an integral domain implies that $A_{V}$ is reduced (\cite{Berkovich}, Proposition 2.3.4). Thus, by \cite{Berkovich}, Proposition 2.1.4, the spectral seminorm on $A_{V}$ is a complete norm, so $v$ is both topologically simple and complete. By Lemma \ref{lem:topsimp:complete}, $A_{V}$ must be a field. By the Nullstellensatz for strictly affinoid algebras, $A_{V}$ is then a finite extension of $K$. Then $\mathcal{H}(v)$ is a finite extension of $K$ (since $A_{V}$ is already a nonarchimedean field, it must be equal to $\mathcal{H}(v)$). Since $A$ injects into $\mathcal{H}(v)$, we conclude that $A$ must be itself a finite extension of $K$.   \end{proof}

\section{Proof of Theorem \ref{Main theorem, introduction}}\label{sec:proof of main theorem}

We can now finally prove Theorem~\ref{Main theorem, introduction} from the introduction. We deduce it from the following more precise statement in characteristic $p$. Recall that we denote by $A^{\coperf}$ the perfect closure of a ring $A$ of characteristic $p$; in particular, for a perfectoid field of characteristic $p$, we denote by $(T_{n, K})^{\coperf}$ the dense subring of the perfectoid Tate algebra $T_{n, K}^{\perfd}$ consisting of $p$-power roots of elements of the usual, Noetherian Tate algebra $T_{n, K}=K\langle X_1,\dots, X_n\rangle$. 
\begin{thm}\label{General statement}Let $K$ be a perfectoid field of characteristic $p$ and let $\mathfrak{m}$ be a maximal ideal in the perfectoid Tate algebra \begin{equation*}T_{n, K}^{\perfd}=K\langle X_{1}^{1/p^{\infty}},\dots, X_{n}^{1/p^{\infty}}\rangle\end{equation*}in $n$ variables, for some $n\geq1$. Set\begin{equation*}L=T_{n, K}^{\perfd}/\mathfrak{m}.\end{equation*}Then there exists a non-negative integer \begin{equation*}l\leq \min(n-\height(\mathfrak{m}\cap (T_{n, K})^{\coperf}), n-1)\end{equation*}and a polyradius $(r_1,\dots, r_{l})\in (0, 1]^{l}$ of length $l$ such that $L$ is a semi-immediate extension of the nonarchimedean field $K_{r_{1},\dots, r_{l}}$. 

Moreover, if \begin{equation*}\mathfrak{m}\cap (T_{n, K})^{\coperf}\neq 0,\end{equation*}then at least one of the radii $r_1,\dots, r_l$ has to be irrational, i.e., $r_{i}\not\in \sqrt{\vert K^{\times}\vert}$. \end{thm}
Before we embark on the proof of Theorem~\ref{General statement}, we establish a simple lemma.
\begin{lemma}
\label{Topologically simple norms and perfection}
Let $A$ be a ring of characteristic $p$, with perfect closure $A^{\coperf}$. 

Then a power-multiplicative norm $\lVert\cdot\rVert$ on $A^{\coperf}$ is topologically simple if and only if its restriction to $A$ is topologically simple.\end{lemma}
\begin{proof}
    Suppose that the restriction is not topologically simple. 
    By Lemma~\ref{Every closed ideal is contained in a closed prime ideal}, there exists a non-zero spectrally reduced prime ideal $\mathfrak{p}$ of the normed ring $(A, \lVert\cdot\rVert\vert_{A})$. 
    By definition of a spectrally reduced prime ideal, there exists a power-multiplicative seminorm $\phi: A\to \mathbb{R}_{\geq0}$ bounded above by $\lVert\cdot\rVert\vert_{A}$ such that $\mathfrak{p}=\ker(\phi)$. 
    We can then extend $\phi$ to a power-multiplicative seminorm $\phi^{1/p^{\infty}}$ on
    $A^{\coperf}$ by setting \begin{equation*}\phi^{1/p^{\infty}}(f^{1/p^{n}})=\phi(f)^{1/p^{n}}\end{equation*}for any $f\in A$, $n\geq0$. 
    By construction, $\phi^{1/p^{\infty}}$ is bounded above by $\lVert\cdot\rVert$. 
    But the kernel of $\phi^{1/p^{\infty}}$ is the prime ideal 
    \begin{equation*}
        \mathfrak{p}^{1/p^{\infty}}=\{\, f\in A^{\coperf}\mid f^{p^{n}}\in \mathfrak{p}~\textrm{for some}~n\geq0\,\},
    \end{equation*}
    which is non-zero since $\mathfrak{p}$ was non-zero. 
    This shows that the normed ring $(A^{\coperf}, \lVert\cdot\rVert)$ has a non-zero closed prime ideal, i..e, it is not topologically simple.
\end{proof}
\begin{proof}[Proof of Theorem \ref{General statement}]We begin by reducing the desired statement to
    the study of a topologically simple continuous absolute value on an affinoid algebra in the
    sense of Tate. For our maximal ideal $\mathfrak{m}\subsetneq T_{n, K}^{\perfd}$, we set \begin{equation*}A_{\infty}\colonequals (T_{n, K})^{\coperf}/(\mathfrak{m}\cap (T_{n, K})^{\coperf}).\end{equation*} 
\begin{claim}\label{claim:null:topsimp}The absolute value on \(L\) induces a topologically simple
        absolute value $v$ on \(A_{\infty}\) which is bounded with respect to the (quotient) norm on
        \(A_{\infty}\) and thus defines an element in the Gelfand spectrum of \(A_{\infty}\).
    \end{claim}

    \begin{proof}[Proof of Claim \ref{claim:null:topsimp}]
        We have the diagram
        \[
        \begin{tikzcd}
        T_{n, K}^{\perfd} \arrow{r}{} \arrow{d}[swap]{} &
        L \arrow{d}{} \\
        (T_{n, K})^{\coperf} \arrow{r}{} &
        A_{\infty}
        \end{tikzcd}
        \]
        of rings.
        Note that the top and left hand maps are continuous.
        Thus, if we endow \(A_{\infty}\) with the topology via
        restriction from \(L\), the diagram is commutative in 
        the category of topological rings.
        
        Now, \(A_{\infty}\)
        is a dense
        subring of $L$, being the image of the dense subring $(T_{n, K})^{\coperf}$ of $T_{n, K}^{\perfd}$ under the continuous surjective map \begin{equation*}T_{n, K}^{\perfd}\to L.\end{equation*}
        In particular, the restriction to $A_{\infty}$ of the absolute value on $L$ is a topologically simple absolute value $v$ on $A_{\infty}$ by Lemma \ref{lem:topsimp:complete}. 
        
        Since \(L\) is a nonarchimedean field, and \(T_{n, K}^{\perfd} \surj L\) is bounded,
        \(v\) pulls back to some point in the Gelfand spectrum of the perfectoid Tate algebra.
        Since the Gelfand spectrum does not change under completion,
        we deduce that the pullback of \(v\) is bounded above by the norm on \((T_{n, K})^{\coperf}\).
        Thus, $v$ is bounded above by the quotient norm on \(A_{\infty}\).
    \end{proof}

    Now, \(\mathfrak{m} \cap (T_{n, K})^{\coperf}\) is the perfect closure of the prime ideal $\mathfrak{m}\cap T_{n, K}$ of $T_{n, K}$ (the perfect closure of $\mathbb{F}_{p}$-algebras inducing a homeomorphism on prime spectra). Thus, $A_{\infty}$ is equal to the perfect closure $A^{\coperf}$ of the affinoid $K$-algebra (in the sense of Tate) 
    \begin{equation*}
        A=T_{n, K}/(\mathfrak{m} \cap T_{n, K}).
    \end{equation*}Since $T_{n, K}$ is excellent (and, in particular, catenary), so is its perfect closure $(T_{n, K})^{\coperf}$. Thus \begin{equation*}n-\height(\mathfrak{m}\cap (T_{n, K})^{\coperf})=\dim(A_{\infty})=\dim(A).\end{equation*}By Lemma~\ref{Topologically simple norms and perfection}, the restriction of $v$ to $A$ is a topologically simple absolute value on $A$, which, by abuse of notation, we again denote by $v$. Let $L'$ be the completion of $A$ with respect to $v$. By our setup, the perfect closure of $L'$ is a dense subfield of $L$, so that $L$ is a semi-immediate extension of $L'$. Thus it suffices to prove that $L'$ is a semi-immediate extension of $K_{r_{1},\dots, r_{l}}$, where $(r_{1},\dots, r_{l})\in (0, 1]^{l}$, where the length $l$ of the polyradius $(r_1,\dots, r_{l})$ satisfies $l\leq \min(\dim(A), n-1)$ and where in the case $l=\dim(A)$ at least one of the radii $r_1,\dots, r_l$ does not belong to the divisible closure of the value group of $K$. To this end, we separate our argument into cases. Note that the Krull dimension of $A$ is always the same as the Krull dimension of $A_{\infty}$ since perfect closure induces a homeomorphism on prime spectra. 
 
\textbf{Case 1: \(A\) has Krull dimension $0$.}

By the Nullstellensatz for affinoid algebras in the sense of Tate, $A$ must be a finite
extension of \(K\), hence complete. Thus $L'=A$ and then $L=L'$ is a finite extension of $K$ since, in this case, $L'$ is already both perfect and complete. 

\textbf{Case 2: \(A\) has Krull dimension $d>0$.}

Let $K\langle T_1,\dots, T_{d}\rangle\hookrightarrow A$ be a Noether normalization of $A$ and let $w$ be the restriction of $v$ to $K\langle T_1,\dots, T_d\rangle$. Since $K\langle T_1,\dots, T_d\rangle \hookrightarrow A$ is finite, the field extension $L'/\mathcal{H}(w)=\mathcal{H}(v)/\mathcal{H}(w)$ is finite, by Lemma~\ref{Finite field extensions and completion}. Hence it suffices to prove that $\mathcal{H}(w)$ is a semi-immediate extension of $K_{r_1,\dots, r_l}$, where either $l<d$ or $l=d$ and at least one of the radii $r_{1},\dots, r_{l}$ is irrational (i.e., does not belong to the divisible closure of the value group of $K$). By Theorem~\ref{Topologically simple rules out Abhyankar}, the restriction $w$ of $v$ to $K\langle T_1,\dots, T_{d}\rangle$ cannot be a rational Abhyankar valuation. If $w$ is an irrational Abhyankar valuation, then, by Lemma~\ref{Temkin's remark and polyradii 4}, $\mathcal{H}(w)=K_{r_1,\dots, r_d}$ and at least one of the radii $r_{1},\dots, r_{d}$ is irrational. If $w$ is not an Abhyankar valuation, then the finitely generated valued extension $K(T_1,\dots, T_d)/K$, with the rank-$1$ valuation induced by $w$, is not Abhyankar and then Lemma~\ref{Temkin's remark and polyradii 3} yields an algebraically independent subset $B=\{x_1,\dots, x_l\}$ of $K(T_1,\dots, T_d)$ of cardinality $l<d$ such that the restricted valuation $w_{B}$ on $K(B)$ is Abhyankar, such that the completion of $K(B)$ with respect to $w_{B}$ is of the form $K_{r_1,\dots, r_l}$ for some polyradius $(r_1,\dots, r_l)\in (0, 1]^{l}$ and such that $\mathcal{H}(w)/\mathcal{H}(w_{B})$ is a semi-immediate extension. 

This concludes the proof of Case 2, and thus also the proof of the theorem.   
\end{proof}

We also include the following restriction that 
rules out (among other things)
the completed residue field of a type II point. 
The argument follows \cite{Gleason22N} precisely.

\begin{prop}
    [\cite{Gleason22N}, Theorem 3.1]
    \label{prop:residue:countable}
    Let \(L\) be a nonachimedean field extension
    of a perfectoid field \(K\) (of characteristic \(p\)).

    Suppose that \(L\)
    is the quotient
    of a perfectoid Tate algebra
    \(K\langle x_{1}^{1 / p^{\infty}},
    \ldots, x_{n}^{1 / p^{\infty}}\rangle \)
    via some surjection \(\varphi\).
    Then \(\widetilde{L}\) has
    a countable basis over \(\widetilde{K}\).
\end{prop}

\begin{proof}
    The proof is exactly the proof in 
    the end of
    \cite[Theorem~3.1]{Gleason22N},
    which we include for completeness.

    Let \(\Gamma\) be the value group of \(L\).
    Then we can embed \(L\) 
    into its field of Mal'cev-Neumann series
    \(\overline{\widetilde{L}}((x^{\Gamma}))\).
    The map 
    \(\varphi\)
    is determined by the \(n\) images
    of the variables \(x_{1}, \ldots, x_{n}\).
    We may write these images as
    \[
    \sum_{j \in \Gamma}^{\infty} a_{i,j}x^{j} 
    \] 
    for \(i \in 1, \ldots, n\) and
    \(a_{i,j} \in \overline{\widetilde{L}}\).
    
    Let \(R^{\prime} = \varphi^{-1}(L^{\circ})\).
    Thus, we have a map 
    \(\tilde{\varphi} \colon
    R^{\prime} \xrightarrow{} \overline{\widetilde{L}}\) 
    via composition with the map
    \(\sum_{}^{} b_{j}x^{q_{j}} \mapsto b_{0}\).
    Let 
    \(S = \{a_{i,j}\}_{i = 1, \ldots, n, j \in \Gamma}\).
    Note that \(S\) is countable, since
    the sets \(\{a_{i,j}\}_{i = k, j \in \Gamma} \) 
    are well-ordered since they 
    are coefficients of a Mal'cev-Neumann series.
    Now, \(\Image(\tilde{\varphi}) \in
    \widetilde{K}[S^{1 / p^{\infty}}] \cap \widetilde{L}
    \ins \tilde{L}\).
    Thus, if \(\widetilde{L}\) has an uncountable
    basis over \(\widetilde{K}\), then 
    \(\overline{\varphi}\) is not surjective,
    and so \(\varphi\) is not surjective.
\end{proof}

\begin{cor}
    \label{cor:residue:nouncountable}
    Let \(L\) be a nonarchimedean field extension of a perfectoid
    field \(K\) such that the residue field
    \(\widetilde{L}\) has positive transcendence degree over \(\widetilde{K}\).
    Assume further that \(\tilde{K}\) is uncountable.
    Then \(L\) is not the quotient field of a perfectoid
    Tate algebra (in any number of variables) over $K$.
\end{cor}

\begin{proof}
    Let \(x \in L\) be a transcendental element.
    Then, the set
    \(\frac{1}{x-c}\) as \(c\) ranges over \(\widetilde{K}\) 
    forms an uncountable \(\widetilde{K}\)-linearly
    independent set in \(\widetilde{L}\),
    and we conclude by Proposition 
    \ref{prop:residue:countable}.
\end{proof}

Using the tilting equivalence, we deduce from Theorem \ref{General statement} our main result describing quotient fields of the perfectoid Tate algebra in arbitrary characteristic.
\begin{thm}\label{Main theorem, in the body of the paper}
    Let $\mathfrak{m}\subsetneq T_{n, K}^{\perfd}$ be a maximal ideal of the perfectoid Tate algebra $T_{n, K}^{\perfd}=K\langle X_{1}^{1 / p^{\infty}}, \ldots, X_{n}^{1 / p^{\infty}}\rangle$ and let \begin{equation*}L=T_{n, K}^{\perfd}/\mathfrak{m}\end{equation*}be the corresponding quotient field. Then there exists a polyradius $(r_1,\dots, r_l)\in (0, 1]^{l}$ such that $L$ is a semi-immediate extension of $K_{r_{1},\dots, r_{l}}$, where \begin{equation*}l\leq \min(n-\height(\mathfrak{m}^{\flat}\cap (T_{n, K^{\flat}})^{\coperf}), n-1),\end{equation*}where $\mathfrak{m}^{\flat}$ is the tilt of $\mathfrak{m}$ in the sense of \cite{Dine22}, \S4. Moreover, if $K$ has uncountable residue field, then none of the radii $r_1,\dots, r_l$ above belongs to the divisible closure of the value group of $K$, and if the residue field is arbitrary, but \begin{equation*}\mathfrak{m}^{\flat}\cap (T_{n, K^{\flat}})^{\coperf}\neq0,\end{equation*}then there exists at least one $i$ such that $r_{i}\not\in \sqrt{\vert K^{\times}\vert}$.
\end{thm}
\begin{proof}If $\mathfrak{m}$ is a maximal ideal of $T_{n, K}^{\perfd}$, then $L=T_{n, K}^{\perfd}/\mathfrak{m}$ is a perfectoid field, by \cite{Kedlaya17}, Corollary 2.9.14, and the tilt $\mathfrak{m}^{\flat}$ of $\mathfrak{m}$ (in the sense of \cite{Dine25}, \S4) is a maximal ideal of $T_{n, K^{\flat}}^{\perfd}$ satisfying \begin{equation*}T_{n, K^{\flat}}^{\perfd}/\mathfrak{m}^{\flat}=L^{\flat},\end{equation*}by \cite{Dine22}, Lemma 4.8. The tilting correspondence for perfectoid fields takes $K_{r_{1},\dots, r_{l}}^{\perfd}$, for any $l\geq0$, to $K_{r_{1},\dots, r_{l}}^{\perfd\flat}$. It also preserves value groups and residue fields, so a perfectoid field over $K$ is a semi-immediate extension of $K_{r_{1},\dots, r_{l}}^{\perfd}$ if and only if its tilt is a semi-immediate extension of $K_{r_{1},\dots, l}^{\flat \perfd}$. We also note that a perfectoid field is a semi-immediate extension of $K_{r_{1},\dots, r_{l}}$ if and only if it is a semi-immediate extension of $K_{r_{1},\dots, r_{l}}^{\perfd}$. Consequently, the theorem follows from Theorem \ref{General statement} and Corollary \ref{cor:residue:nouncountable}. \end{proof}
\begin{cor}\label{Corollary of the main theorem, in the body of the paper}Let $K$ be a perfectoid field and let $n\geq1$ be an integer. No nonarchimedean field extension $L$ of $K$ which is a semi-immediate extension of a nonarchimedean field of the form $K_{r_{1},\dots, r_{n}}$, or which admits a nonarchimedean field of the form $K_{r_{1},\dots, r_{n}}$ as a semi-immediate extension, can be a quotient field of the perfectoid Tate algebra in $n$ variables.\end{cor}
\begin{proof}This follows from the bound $l\leq n-1$ in Theorem \ref{Main theorem, in the body of
    the paper} by comparing $d_{K}(L)$ with the analogous quantity for a nonarchimedean field extension of the form $K_{r_{1},\dots, r_{l}}/K$ with $l\leq n-1$.
\end{proof}
\begin{cor}\label{Two variables and large value group}Let $K$ be a perfectoid field with $\vert K^{\times}\vert=\mathbb{R}_{>0}$ and let \begin{equation*}\mathfrak{m}\subsetneq T_{2, K}^{\perfd}\end{equation*}be a maximal ideal of the perfectoid Tate algebra in two variables over $K$. If \begin{equation*}\mathfrak{m}^{\flat}\cap (T_{n, K^{\flat}})^{\coperf}\neq0,\end{equation*}then the quotient field $T_{2, K}^{\perfd}/\mathfrak{m}$ is a semi-immediate extension of $K$. In particular, if $K$ is also algebraically closed and spherically complete, then $T_{n, K}^{\perfd}/\mathfrak{m}=K$. \end{cor}
\begin{proof}By Theorem \ref{Main theorem, in the body of the paper}, $L=T_{2, K}^{\perfd}/\mathfrak{m}$ is either a semi-immediate extension of $K$ or a semi-immediate extension of $K_{r}$ for some $r$ which does not belong to the divisible closure of the value group of $K$. However, the value group of $K$ is $\mathbb{R}_{>0}$, so the second possibility cannot occur.\end{proof}
More generally, if the value group of $K$ is $\mathbb{R}_{>0}$, the proof of Theorem \ref{Main theorem, in the body of the paper} yields the following sharpening of the upper bound in that theorem.
\begin{prop}\label{Any number of variables and large value group}Let $K$ be a perfectoid field with $\vert K^{\times}\vert=\mathbb{R}_{>0}$, let $n\geq1$ and let \begin{equation*}\mathfrak{m}\subsetneq T_{n, K}^{\perfd}\end{equation*}be a maximal ideal of the perfectoid Tate algebra in $n$ variables over $K$. Then there exists a polyradius $(r_1,\dots, r_l)\in (0, 1]^{l}$ of length \begin{equation*}l\leq n-\height(\mathfrak{m}^{\flat}\cap (T_{n, K^{\flat}})^{\coperf})-1\end{equation*}such that \begin{equation*}L=T_{n, K}^{\perfd}/\mathfrak{m}\end{equation*}is a semi-immediate extension of $K_{r_1,\dots, r_l}^{\perfd}$. \end{prop}
\begin{proof}By tilting, the assertion reduces once again to the case when $K$ is of characteristic $p$. Since the value group of $K$ is all of $\mathbb{R}_{>0}$, every Abhyankar valuation on any affinoid algebra $A$ over $K$ must be rational. Thus, by Theorem \ref{Topologically simple rules out Abhyankar}, no topologically simple continuous absolute value on any strictly $K$-affinoid algebra $A$ can be Abhyankar. As in the proof of Theorem \ref{General statement}, the absolute value on $L$ restricts to a continuous absolute value $v$ (over $K$) on the strictly affinoid algebra $A$ given by the image of $T_{n, K}$ in $L$, and $L$ is a semi-immediate extension of $\mathcal{H}(v)$. Let $K\langle T_1,\dots, T_d\rangle \hookrightarrow A$ be a Noether normalization of $A$. Since $\mathfrak{m}\cap (T_{n, K})^{\coperf}\neq0$, we have $d<n$. Since the topologically simple absolute value $v$ on $A$ cannot be Abhyankar, its restriction to $K\langle T_1,\dots, T_d\rangle $ is not Abhyankar, by Lemma \ref{Finite extensions and Abhyankar valuations}. By Lemma~\ref{Temkin's remark and polyradii 3} applied to the purely transcendental, finitely generated valued field extension $K(T_1,\dots, T_d)/K$, this means that there exists an algebraically independent subset $B=\{x_1,\dots, x_l\}$ of $K(T_1,\dots, T_d)$ of cardinality $l<d$ such that the completion of $K(B)$ with respect to the restriction of the valuation on $\mathcal{H}(v)$ is of the form $K_{r_1,\dots, r_l}$ for some polyradius $(r_1,\dots, r_l)\in (0, 1]^{l}$ of length $l$ and such that $\mathcal{H}(v)$ is a semi-immediate extension of $K_{r_1,\dots, r_l}$. Since $L$ was a semi-immediate extension of $\mathcal{H}(v)$, this proves that $L$ is a semi-immediate extension of $K_{r_1,\dots, r_l}^{\perfd}$ for $l<d=n-\height(\mathfrak{m}\cap (T_{n, K})^{\coperf})$, as claimed. \end{proof} 

\section{On Gleason's argument}\label{sec:Gleason}

In this section, we give a slightly simpler proof of Gleason's argument.
Then, we generalize this argument to show

\begin{thm}
    \label{thm:multivar:gleason}
    Let $n\geq1$ and let \((r_{1}, \ldots, r_{l}) \in (0,1]^{l}\) be a free polyradius of length $l\leq n$. 
    The field \(K_{r_{1}, \dots, r_{l}}^{\perfd}\)
    is a quotient field of 
    \(K\langle T_{1}^{1 / p^{\infty}}, \ldots, T_{n+2}^{1 / p^{\infty}}\rangle\).
    Furthermore, there exists a field \(K^{\prime}\) 
    for which \(K_{r_{1}, \ldots, r_{n}}^{\perfd}\) is a
    finite separable extension
    such that \(K^{\prime}\) 
    is a quotient field of 
    \(K\langle T_{1}^{1 / p^{\infty}}, \ldots, T_{n+1}^{1 / p^{\infty}} \rangle \).  
\end{thm}
We begin by recalling 
the notion of being \((q,s)\)-adapted from \cite{Gleason22N}.
We adapt the definition to work in \(K_{r_{1}, \ldots, r_{n}}\), for a free polyradius $(r_1,\dots, r_n)$, 
rather than just \(K_{r}\) for a single radius $r\not\in\sqrt{\vert K^{\times}\vert}$.
Our convention is that \(K_{r_{1}, \ldots, r_{n}}\) 
(and thus \(K_{r_{1}, \ldots, r_{n}}^{\perfd}\))
is generated (topologically) by \(x_{1}, \ldots, x_{n}\),
and for \(q = (q^{1}, \ldots, q^{n})\),
\(x^{q} \colonequals x_{1}^{q^{1}}\cdot\ldots \cdot x_{n}^{q^{n}}\).
Furthermore,
\(r \colonequals (r_{1}, \ldots, r_{n})\),
and \(r^{q} \colonequals r_{1}^{q^{1}}\cdot \ldots \cdot r_{n}^{q^{n}}\).
Finally, let \(|\cdot|_{r}\) 
denote the \((r_{1}, \ldots, r_{n})\)-Gauss norm.
This choice of convention makes the notation compatible
with the one-variable case; 
that is, if \(n = 1\), then one has \(q \in \mathbb{Z}[\frac{1}{p}]\),
\(r \in (0,1]\), and e.g. \(r^{q}\) is exponentiation
of real numbers.
Thus, the reader may safely assume that \(n = 1\)
for a first reading and the proof still holds.

\begin{defn}
    Let \(0 < s < 1\) and \(q \in \mathbb{Z}[1 / p]^{n}\).
    An element \(\beta = \sum_{j \in \mathbb{Z}[1 / p]^{n}}^{} b_{j}x^{j} 
    \in K_{r_{1}, \ldots, r_{n}}^{\perfd}\) is 
    \textit{\((q,s)\)-adapted} if
    \begin{enumerate}[(1)]
        \item \(s < |\beta|_{r} \leq 1\)
        \item  \(\argnorm \beta = q\), that is,  
            \(|b_{q}|r^{q} = |\beta|_{r}\).
        \item \(|\beta - b_{q}x^{q}|_{r} \leq s|\varpi|\).
    \end{enumerate}
\end{defn}

\begin{figure}[h]
    \label{fig:qsadapted:drawing}
\centering
\includegraphics[width=0.75\textwidth]{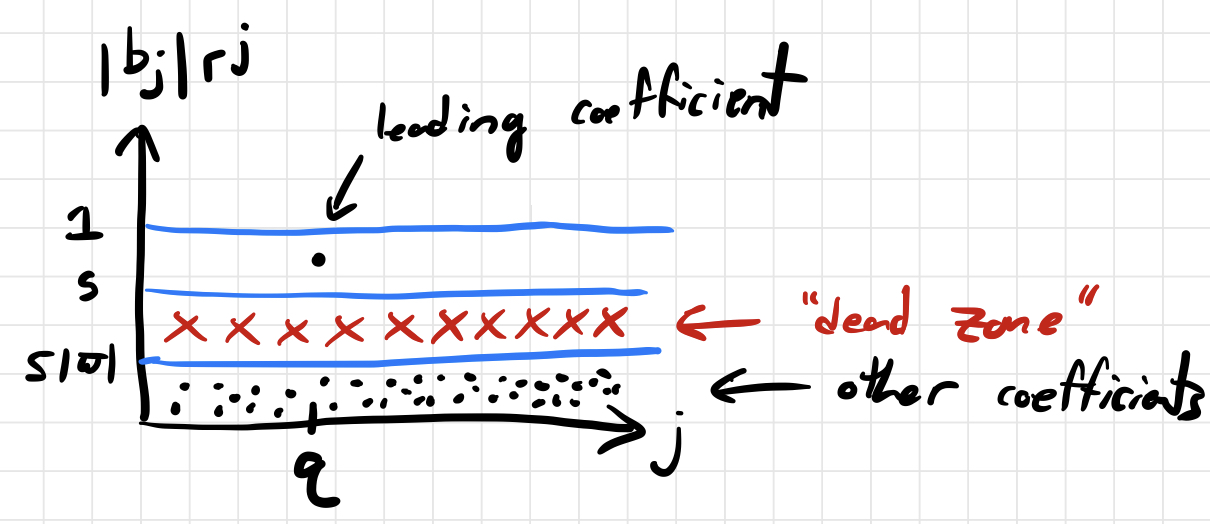}
\caption{A picture of \((q,s)\)-adapted
    when \(n=1\)
    for the element \(\beta = \sum_{j \in \mathbb{Z}[1 / p]}^{} b_{j}x^{j} \in K_{r}^{\perfd}\).
    Note that the picture is lying slightly, the \(j\)-axis should extend in the
    negative direction as well.
}
\end{figure}

The notion of being \((q,s)\)-adapted may be drawn pictorially in Figure 
\ref{fig:qsadapted:drawing}.

\subsection{Division Algorithm}

\begin{defn}
    Loosely following the notation of \cite[Section~2.1]{garzella-2024-perfectoid-uncountable},
    for \(\beta = \sum_{q \in \mathbb{Z}[\frac{1}{p}]^{n}}^{} b_{q}x^{q}
    \in K_{r_{1}, \ldots, r_{n}}^{\perfd}\),
    define \(\res^{M \leq} \beta\) by
    \[
        \res^{M \leq} \beta = \sum_{q \in \mathbb{Z}[\frac{1}{p}]^{n}, 
        M \leq |b_{q}| r^{q}}^{} b_{q}x^{q}
    \] 
\end{defn}

In the previous definition, ``res'' stands for ``restriction''.
Note that as \(\beta \in K_{r_{1}, \ldots, r_{n}}^{\perfd}\),
\(\res^{M \leq} \beta\) is always a finite sum.

\begin{figure}[h]
    \label{fig:res_greater:drawing}
\centering
\includegraphics[width=0.75\textwidth]{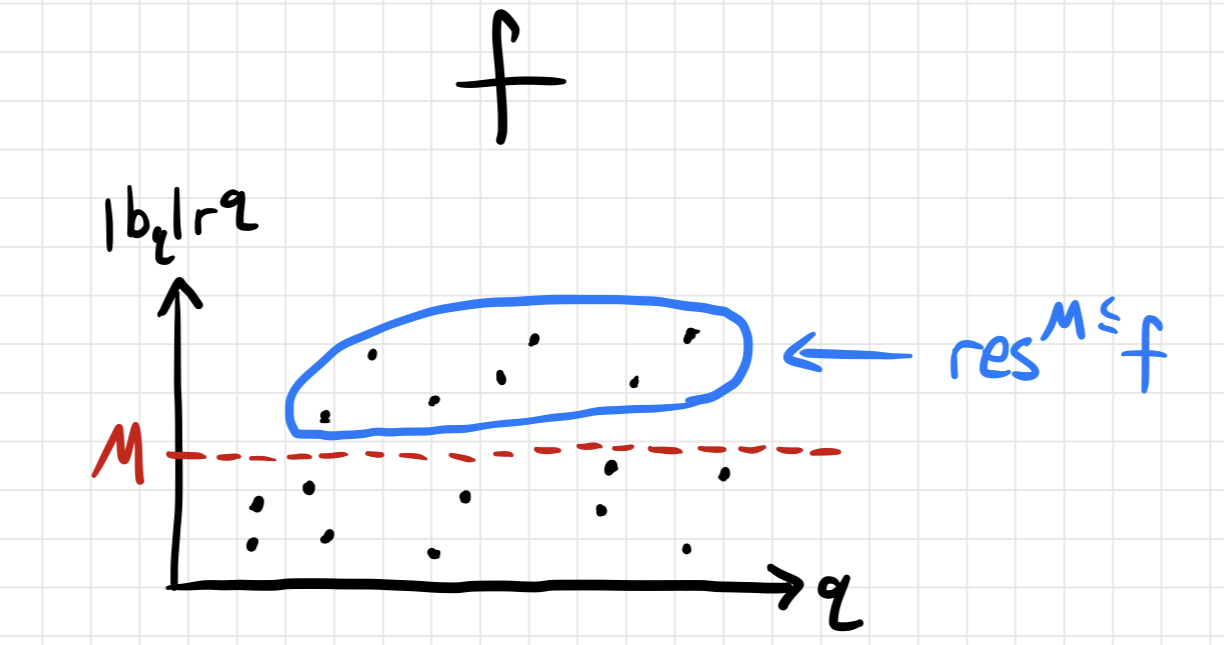}
\caption{A picture of \(\res^{M \leq} f\)
    when \(n=1\),
    which consists of the (finitely many) terms of \(f\) 
    who have norm greater than \(M\).
}
\end{figure}

\begin{prop}
    \label{prop:qsadapted:surj}
    [\cite{Gleason22N}, Lemma 4.3]
    Let \(R\) be a complete Tate \(K\)-algebra and let 
    \(\phi \colon R \xrightarrow{} K_{r_{1}, \ldots, r_{n}}^{\perfd}\) be a 
    continuous map.
    Suppose that there exists \(0 < s < 1\) such that, for all $q$,
    there is an element \(f \in R^{\circ}\) such that
    \(\phi(f)\) is \((q,s)\)-adapted.
    Then \(\phi\) is surjective.
\end{prop}

The proof of \cite[Lemma~4.3]{Gleason22N} applies in this context,
as that proof never uses any properties about the source
other than complete Tate. 
For completeness we give the argument.

\begin{lem}
    [Division algorithm for \(K_{r}^{\perfd}\)]
    \label{lem:res:divisionalg}
    Let \(\phi\) be as in Proposition \ref{prop:qsadapted:surj}. 
    Then for any \(\beta \in K_{r}^{\perfd}\),
    such that \(|\varpi|^{m+1}s \leq |\beta|_{r} \leq |\varpi^{m}|s\)
    for \(0 \leq m\),
    there exists an \(f \in R^{\circ}\) such that 
    \(|\phi(f) - \beta|_{r} \leq |\varpi|^{m+1}s\).
    Moreover, \(|f| \leq |\varpi|^{m}\).
    In other words, 
    there exists $f\in R^{\circ}$ such that 
    \(\res^{|\varpi|^{m+1}s \leq} \phi(f)
    = \res^{|\varpi|^{m+1} s \leq} \beta\).
\end{lem}

\begin{proof}
   We have that \(\res^{|\varpi|^{m+1} s \leq} \beta = \sum_{i=1}^{n} b_{i}x^{q_{i}}\)
   for some \(q_{i} \in \mathbb{Z}[\frac{1}{p}]\).
   Take elements \(f_{i}\) 
   such that \(\phi(f_{i})\) is \((q_{i}, s)\)-adapted.
   Then \(\phi(f_{i}) = c_{i}x^{q_{i}} + f_{i}^{\prime}\),
   for some \(c_{i}\) such that \(s \leq |c_{i}|r^{q_{i}}\)
   and \(|f_{i}^{\prime}|_{r} \leq |\varpi|s\).
   
   By our assumption on \(\beta\), \(|\frac{b_{i}x^{q_{i}}}{\varpi^{m}}|_{r} \leq s < |c_{i}x^{q_{i}}|_{r}\),
   so \(|\frac{b_{i}}{c_{i}}| < |\varpi|^{m} \).
   Let \(d_{i} = \frac{b_{i}}{c_{i}} \in K\).
   Now take \(f = \sum_{i=1}^{n} d_{i}f_{i}\).
   We see that by the assumption
   on \(\phi\), all elements \(f_{i}\) are in \(R^{\circ}\).
   We get that \(f\) is in \(R^{\circ}\)
   and has norm less than \(|\varpi|^{m}\).
  
   By construction, we have that 
   \[
   \phi(f) = \res^{|\varpi|^{m+1} s \leq} \beta + \sum_{i=1}^{n} d_{i}f_{i}^{\prime} 
   .\] 
   Now, since \(|d_{i}| \leq |\varpi|^{m}\)
   and \(|f_{i}^{\prime}| \leq |\varpi|s\) 
   for all \(i\), 
   we have 
   \(|d_{i}f_{i}^{\prime}| \leq |\varpi|^{m+1}s\) 
   and we conclude.
\end{proof}

\begin{proof}[Proof of Proposition \ref{prop:qsadapted:surj}]
    Let 
    \[\beta = \sum_{q \in \mathbb{Z}[\frac{1}{p}]}^{} b_{q}x^{q} \in K_{r}^{\perfd}.\]
    By replacing \(\beta\) with \(\varpi^{k} \beta\),
    we may assume that \(|\beta|_{r} \leq s\).
    We construct a Cauchy sequence \(f_{m} \in R\)
    and a null sequence \(\beta_{m} \xrightarrow{} 0\) 
    in \(K_{r}^{\perfd}\) such that
    \(\phi(f_{m}) = \beta - \beta_{m}\).

    More precisely, our null sequence \(\beta_{m}\) 
    will have the property that \(|\beta_{m}| \leq |\varpi|^{m}s\).

    Let \(\beta_{0} = \beta\) and \(f_{0} = 0\).
    To construct \(\beta_{m}\) and \(f_{m}\),
    note that
    by Lemma \ref{lem:res:divisionalg},
    there exists an element \(g_{m}\) such that 
    \(|\beta_{m} - \phi(g_{m})|_{r}\leq |\varpi|^{m+1}s\) and $\vert g_{m}\vert\leq \vert\varpi\vert^{m}$.
    Thus we may take \(f_{m+1} = f_{m} + g_{m}\),
    and \(\beta_{m+1} = \beta_{m} - \phi(g_{m})\).
    Now \(|f_{m+1} - f_{m}| \leq |\varpi|^{m}\), so \(f_{m}\) is Cauchy.
    Finally, \(\phi(f_{m+1}) = \phi(f_{m}) + \phi(g_{m})
    = \beta - \beta_{m} + \phi(g_{m}) = \beta - \beta_{m+1}\).
\end{proof}

\subsection{Gleason Elements}

\begin{defn}
    Let \(J \ins \mathbb{Z}[\frac{1}{p}]^{n}\). Let \(\omega\) be a well-order of \(J\).
    A \textit{Gleason element} for \(J\) is an element
    \(\beta \in K_{r}^{\perfd}\) such that 
    for every \(m \in \mathbb{N}\) there exists 
    \begin{itemize}
        \item a monomial \(W_{m} = w_{m}x^{\omega(m)}\) in 
            \(K(x_{1}^{1 / p^{\infty}}, \ldots, x_{n}^{1 / p^{\infty}})\)
        \item an \(\varepsilon_{m} \in K\) and
            \(b_{m} \in \mathbb{N}\) such that 
            \(\varepsilon_{m}\) depends on \(b_{1}, \ldots, b_{m-1}\) 
            and \(b_{m}\) depends on \(\varepsilon_{m}\).
        \item for \(1 \leq i \leq m\), a \(d_{m,i} \in K^{\circ}\) depending on 
            \(W_{i}\), \(\varepsilon_{m}\), and 
            \(b_{m}\)
    \end{itemize}
    such that
    \begin{equation}
        \label{eq:gleason:elt}
        \left(\varepsilon_{m} \beta - \sum_{i=1}^{m} d_{m,i} W_{i}^{p^{b_{i}}}\right)^{1 / p^{b_{m}}}
    \end{equation}
    is \((\omega(m),s)\)-adapted.
    Furthermore, we require that \(b_{m} \xrightarrow{} \infty\) 
    as \(m \xrightarrow{} \infty\)
    and the exponents of
    the monomials \(W_{m}^{p^{b_{m}}}\) are all distinct.
\end{defn}

\begin{figure}[h]
    \label{fig:res_greater:drawing}
\centering
\includegraphics[width=0.75\textwidth]{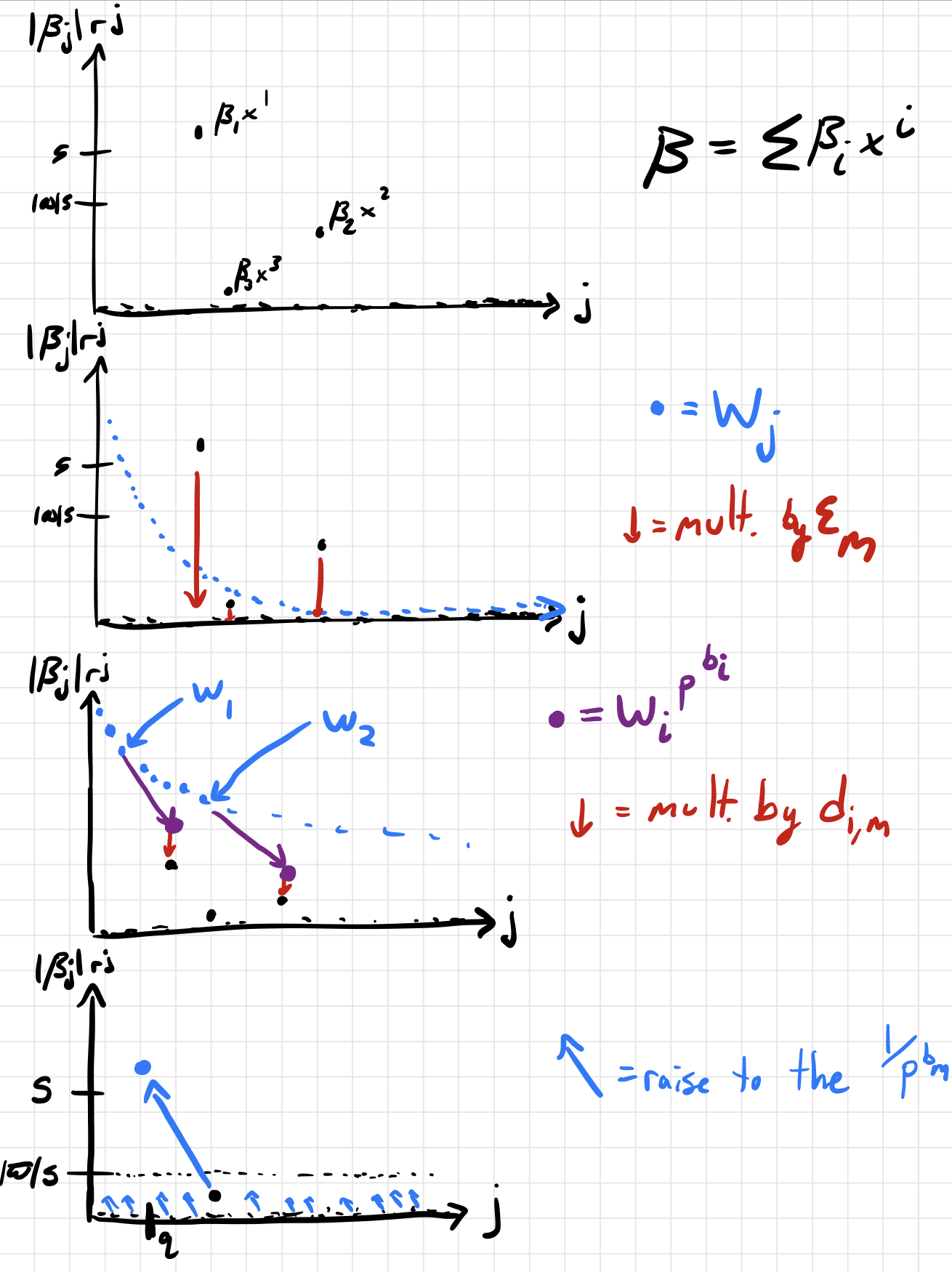}
\caption{A picture of a Gleason element, for \(n=1\) 
    and \(J = \mathbb{Z}[\frac{1}{p}]_{\geq 0}\).
    Some \(p\)-th root of
    a linear combination of \(\beta\) and monomials
    is \((q,s)\)-adapted.
    The \((q,s)\)-adapted elements and the corresponding
    linear combination are obtained
    recursively by eliminating previous terms.
}
\end{figure}

\begin{lem}
    \label{lem:gleason:argnorm}
    Let \(\beta \in K_{r_{1}, \ldots, r_{n}}\) be a Gleason element.
    Then for any \(m\),
    \[
    \varepsilon_{m}\beta
    - \sum_{i=1}^{m} d_{m,i}W_{i}^{p^{b_{i}}} 
    \] 
    attains its norm at \(c_{m}W_{m}^{p^{b_{m}}}\),
    for \(c_{m} =  w_{m}^{p^{b_{m}}}\).
\end{lem}

\begin{proof}
    This follows by the fact 
    that the term in Equation \ref{eq:gleason:elt}
    is \((q,s)\)-adapted.
    In particular, the term
    \(c_{m}W_{m}^{p^{b_{m}}}\) has
    norm greater than \(s^{p^{b_{m}}}\),
    while all other
    terms have norm less than \(|\varpi|^{p^{b_{m}}} s^{p^{b_{m}}}\).
\end{proof}

\begin{lem}
    Let \(\beta\) be a Gleason element.
    Let \(c_{m}\) be as in Lemma \ref{lem:gleason:argnorm}.
    Then \(|c_{m}W_{m}^{p^{b_{m}}}| <
    \left| \sum_{i=1}^{m-1} d_{m,i}W_{i}^{p^{b_{i}}} \right|\)
\end{lem}

\begin{proof}
    By Lemma \ref{lem:gleason:argnorm}
    applied to the case \(m=1\),
    the norm of \(\varepsilon_{m}\beta\) 
    must be achieved at \(c^{\prime} W_{1}^{p^{b_{1}}}\) 
    for some \(c^{\prime} \in K\).
    In fact, \(\varepsilon_{m}s \leq |c^{\prime} W_{1}^{p^{b_{1}}}|\),
    while all other terms of 
    \(\varepsilon_{m}\beta\) have norm below
    \( |\varepsilon_{m}\varpi|s\).
    Thus, if the conclusion is contradicted, we
    have 
    that \(W_{1}^{p^{b_{1}}}\) and \(W_{m}^{p^{b_{m}}}\)
    are multiples of the same monomial,
    contradicting the fact that all such
    monomials have distinct exponents.
\end{proof}

\begin{lem}
    \label{lem:gleason:epsilon}
    Let \(\beta\) be a Gleason element.
    Then for any \(m\), 
    \(s \leq |\varepsilon_{m}^{1 / p^{b_{m}}}|\).
\end{lem}

\begin{proof}
    By Lemma \ref{lem:gleason:argnorm},
    we have that \(s^{p^{b_{m}}} \leq |c_{m}W_{m}^{p^{b_{m}}}|\)
    for some \(c_{m} \in K\).
    Now, we have
    \[
    |\varepsilon_{m}\beta| =
    \left|\beta^{\prime} + c_{m}W_{m}^{p^{b_{m}}} + \sum_{i=1}^{m-1} d_{m,i}W_{i}^{p^{b_{i}}}\right|
    \leq \sup \left\{|\beta^{\prime}|, |c_{m}W_{m}^{p^{b_{m}}}|,
    \left| \sum_{i=1}^{m-1} d_{m,i}W_{i}^{p^{b_{i}}}  \right|\right\},
    \] 
    where \(\beta^{\prime}\) has norm less than 
    \(|\varpi|^{p^{b_{m}}} s^{p^{b_{m}}}\).
    Thus, we can conclude that 
    \(|c_{m}W_{m}^{p^{b_{m}}}| \leq |\varepsilon_{m}\beta|\).
    Finally, WLOG \(|\beta| < 1\) so 
    \(|\varepsilon_{m}\beta| \leq |\varepsilon_{m}|\).
    Thus, 
    \(s^{p^{b_{m}}} \leq |\varepsilon_{m}|\), and we conclude.
\end{proof}

\begin{prop}
    [Structure of Gleason elements]
    Let \(\beta \in K_{r_{1}, \ldots, r_{n}} \) be a Gleason element.
    Then \(\beta = \sum_{i=1}^{\infty} \beta_{i}W_{i}^{p^{b_{i}}}\)
    for some coefficients \(\beta_{i} \in K\).
    Let \(W_{i} = w_{i}x^{q_{i}}\) for \(w_{i} \in K\),
    \(q_{i} \in \mathbb{Z}[\frac{1}{p}]^{n}\).
    Then \(\{q_{i} \colon i \in \mathbb{N}\} = J\),
    i.e. the sequence \(\{q_{i}\}_{i}\) is
    a well-order of \(J\).
\end{prop}

\begin{proof}

    First of all, Equation \ref{eq:gleason:elt}
    shows that 
    \(\varepsilon_{m}\beta - \sum_{i=1}^{m-1} d_{m,i}W_{i}^{p^{b_{i}}} \) 
    has a term which is \(c_{m}W_{m}^{p^{b_{m}}}\)
    The fact that \(\omega \colon J \xrightarrow{} \mathbb{N}\)
    is an injection.
    Then we see that 
    \[
    \beta = \frac{\sum_{i=1}^{m-1} d_{m,i}W_{i}^{p^{b_{i}}} + 
    c_{m}W_{m}^{p^{b_{m}}} + \beta^{\prime p^{b_{m}}}}{\varepsilon_{m}} 
    \] 
    for some \(\beta^{\prime}\) with \(|\beta^{\prime}| \leq |\varpi|s\).
    Thus, \(\frac{\beta^{\prime p^{b_{m}}}}{\varepsilon_{m}}\)
    has norm less than \(\frac{|\varpi|^{p^{b_{m}}}s^{p^{b_{m}}}}{|\varepsilon_{m}|}\).
    By Lemma \ref{lem:gleason:epsilon},
    this must be less than \(|\varpi|^{p^{b_{m}}}\).
    Thus, the sequence
    \[
        \left\{ \frac{1}{\varepsilon_{m}} \sum_{i=1}^{m-1} d_{m,i}W_{i}^{p^{b_{i}}} 
        + \frac{1}{\varepsilon_{m}} c_{m}W_{m}^{p^{b_{m}}}  \right\}_{m}
    \] 
    is Cauchy with limit \(\beta\),
    and moreover each of the exponents 
    \(W_{i}^{p^{b_{i}}}\) is distinct.
    The first claim follows.

    The second claim follows from the first
    and the \((\omega(m), s)\)-adaptedness 
    of Equation \ref{eq:gleason:elt}.
\end{proof}

Our main use for Gleason elements is to ``package''
\((q,s)\)-adapted elements to ensure that they are
in the image of a map.

\begin{prop}
    \label{prop:gleason:surj}
    Let \(R\) be a perfectoid Tate ring.
    Let \(\varphi \colon R \xrightarrow{} K_{r_{1}, \ldots, r_{n}}^{\perfd}\) 
    be a continuous map.
    Assume that there exist Gleason elements 
    \(g_{1}, \ldots, g_{n}\) for 
    sets \(J_{1}, \ldots, J_{n}\)
    such that \(\bigcup_{k}^{} J_{k} = \mathbb{Z}[\frac{1}{p}]^{n}\).
    Furthermore, assume that all the \(W_{i}\) for each \(g_{k}\) 
    are in the image \(\varphi(R^{\circ})\).
    Then \(\varphi\) is surjective.
\end{prop}

\begin{proof}
    The assumption implies that a \((q,s)\)-adapted
    element is in the image of \(\varphi\)
    for every \(q \in \mathbb{Z}[\frac{1}{p}]^{n}\).
    Then we conclude by 
    Proposition~\ref{prop:qsadapted:surj}.
\end{proof}

We now give the most basic construction of a Gleason element.

\begin{prop}
    \label{prop:simple:gleason}
    Let \(J = \mathbb{Z}[\frac{1}{p}]_{\geq 0}\),
    and let \(\omega\) be a well-ordering of \(J\).
    There exists a Gleason element \(G_{+}
    \in K_{r_{1}}^{\perfd}\) 
    for \(J\) which has monomials 
    \(W_{i} = x_{1}^{\omega(i)}\).
\end{prop}

\begin{proof}
    We follow \cite[Proposition~4.4]{Gleason22N}.
    For convenience, since there is only one variable, 
    let \(x \colonequals x_{1}\).
    First, choose for all \(m \in \mathbb{N}\),
    choose \(e_{m} \in K\) such that 
    \(s < |e_{m}x^{\omega(m)}|_{r} < 1\).
    Let \(\alpha_{m} \colonequals e_{m}x^{\omega(m)}\).
    Now, choose \(\varepsilon_{m}\) and \(b_{m}\) 
    (\(\varepsilon_{m}\) depending on
    \(b_{1}, \ldots, b_{m-1}\), \(b_{m}\) 
    depending on \(\varepsilon_{m}\))
    so that

    \begin{enumerate}[(1)]
        \item \(b_{m}\) is large enough that
            \(|\alpha_{m}|^{p^{b_{m}}} < |\varpi|^{m}\)
        \item \(\varepsilon_{m}\)
            is small enough that
            \(|\varepsilon_{m} \alpha_{i}^{p^{b_{i}}}| \leq |x^{\omega(i)}|^{p^{b_{i}}}\)
            for \(1 \leq i < m\).
        \item \(b_{m}\) is big enough that 
            \(s < |\varepsilon_{m}^{1 / p^{b_{m}}} \alpha_{m}| < 1 \)
        \item \(b_{m}\) is large enough that 
            \(|\alpha_{m}^{p^{b_{m} - b_{i}}}| < |\varpi| s\)
            for all \(1 \leq i < m\)
        \item \(b_{m}\) is large enough so that
            \(\omega(m)p^{b_{m}}\) is distinct
            from all \(\omega(i)p^{b_{i}}\) 
            for \(1 \leq i < m\).
    \end{enumerate}

    Finally, we let \(G_{+}\) be \(\sum_{i}^{} \alpha_{i}^{p^{b_{i}}} \).
    Now, we need to show that 
    there exist \(d_{m,i}\) such that
    \(H \colonequals \left(  \varepsilon_{m}G_{+} - \sum_{i=1}^{m-1} d_{m,i}W_{i}^{p^{b_{i}}} \right)^{1 / p^{b_{m}}}\)
    is \((\omega(m),s)\)-adapted.

    Let \(d_{m,i}\) be \(\varepsilon_{m}e_{i}^{p^{b_{i}}}\).
    Then, \(d_{m,i}W_{i}^{p^{b_{i}}} = \varepsilon_{m} e_{m}^{p^{b_{i}}} x^{\omega(i) p^{b_{i}}}\)
    is precisely the \(i\)-th term of \(\varepsilon_{m} G_{+}\).
    By (2), \(|d_{m,i}| < 1\).
    Now, by (3) the term 
    \(\varepsilon_{m}^{1 / p^{b_{m}}} e_{m}x^{\omega(m)}\) 
    has norm between \(s\) and \(1\),
    while all other terms of the expression \(H\) 
    have norm less than \(|\varpi|s\) by (4).
    Finally, condition (5) ensures that \(G_{+}\) satisfies
    the unique exponents property in the definition of Gleason element.
    Thus, \(H\) is \((\omega(m),s)\)-adapted.
\end{proof}

The construction in this case readily applies
to the most general case:

\begin{prop}
    \label{prop:gleason:multivar}
    Let \(V_{1}, \ldots, V_{\ell}\) be monomials
    in \(K(x_{1}^{1 / p^{\infty}}, \ldots, x_{n}^{1 / p^{\infty}})\)
    such that \(|V_{k}|_{r} < 1\),
    and the exponents of \(V_{1}, \ldots, V_{n}\) 
    are \(\mathbb{Z}[\frac{1}{p}]_{\geq 0}\)-linearly 
    independent.
    Let \(J\) be the \(\mathbb{Z}[\frac{1}{p}]\)-lattice
    spanned by  the \(V_{k}\), \(k = 1, \ldots, n\).
    Then there exists a Gleason element 
    in \(K_{r_{1}, \ldots, r_{n}}^{\perfd}\)
    for \(J\).
\end{prop}

\begin{proof}
    Let \(W_{i}\) 
    be constructed as follows:
    Let \(v_{k} \in \mathbb{Z}[\frac{1}{p}]^{n}\) be the 
    exponent tuple of \(V_{k}\).
    Write \(\omega(i) = \sum_{k=1}^{\ell} h_{k}v_{k}\) as a sum of 
    the \(v_{k}\), which by assumption
    is unique.
    Then \(h_{k}\) are in \(\mathbb{Z}[\frac{1}{p}]_{\geq 0}\).
    Then, let \(W_{i} = V_{1}^{h_{1}}\dots V_{\ell}^{h_{\ell}}\).
    Note that \(|W_{i}| \leq 1\).
    Now, we mimic the construction of Proposition \ref{prop:simple:gleason}.

    First, let \(\alpha_{m} = e_{m}W_{m}\) for all \(m\)
    such that \(s \leq |e_{m} W_{m}| < 1\).
    Now, choose \(\varepsilon_{m}\) and \(b_{m}\) 
    (\(\varepsilon_{m}\) depending on
    \(b_{1}, \ldots, b_{m-1}\), \(b_{m}\) 
    depending on \(\varepsilon_{m}\))
    so that

    \begin{enumerate}[(1)]
        \item \(b_{m}\) is large enough that
            \(|\alpha_{m}|^{p^{b_{m}}} < |\varpi|^{m}\)
        \item \(\varepsilon_{m}\)
            is small enough that
            \(|\varepsilon_{m} \alpha_{i}^{p^{b_{i}}}| \leq |W_{i}|^{p^{b_{i}}}\)
            for \(1 \leq i < m\).
        \item \(b_{m}\) is big enough that 
            \(s < |\varepsilon_{m}^{1 / p^{b_{m}}} \alpha_{m}| < 1 \)
        \item \(b_{m}\) is large enough that 
            \(|\alpha_{m}^{p^{b_{m} - b_{i}}}| < |\varpi| s\)
            for all \(1 \leq i < m\)
        \item \(b_{m}\) is large enough so that
            \(\omega(m)p^{b_{m}}\) is distinct
            from all \(\omega(i)p^{b_{i}}\) 
            for \(1 \leq i < m\).
    \end{enumerate}

    As before, we define \(\beta = \sum_{i}^{} \alpha_{i}^{p^{b_{i}}}\).
    Also as before, (1) ensures that \(\beta\) is convergent,
    while (2)-(5) say that \(\beta\) is a Gleason element,
    exactly as above.
\end{proof}

We also have 

\begin{prop}
    \label{prop:gleason:minus}
    Let \(c \in K\) be such that 
    \(s < \vert cx^{-1}\vert_{r} < 1\). Let \(J = \mathbb{Z}[\frac{1}{p}]_{\leq 0}\).
    There exists a Gleason element \(G_{-}
    \in K_{r_{1}}\)
    such that \(G_{-}\) 
    is in the image of
    the map 
    \begin{align*}
        \phi \colon K\langle T_{1}^{1 / p^{\infty}}\rangle  &\longrightarrow K_{r_{1}}^{\perfd} \\
        T_{1} &\longmapsto cx^{-1}
    .\end{align*}
\end{prop}

\begin{proof}
    By Proposition \ref{prop:gleason:multivar},
    there exists a Gleason element 
    \(G_{-}\) for \(J\).
    Fix \(\epsilon > 0\) 
    such that \(\epsilon < |cx^{-1}| - s\).
    In the proof of Proposition \ref{prop:gleason:multivar}
    choose \(e_{m}\) such that 
    \(s < |e_{m}| < |c x^{-1}| - \epsilon\).
    Now, the element 
    \[
    \sum_{m=1}^{\infty} e_{m}T_{1}^{\omega(m) p^{b_{m}}} 
    \] 
    maps to \(G_{-}\), and we are done.
\end{proof}

We recover the following result, which is essentially
the same as Gleason's:

\begin{cor}
    The map
    \begin{align*}
        \varphi : R \colonequals K\langle T_{1}^{1 / p^{\infty}}, T_{2}^{1 / p^{\infty}}, T_{3}^{1 / p^{\infty}}\rangle &\longrightarrow K_{r_{1}} \\
        T_{1} &\longmapsto x \\
        T_{2} &\longmapsto cx^{-1} \\
        T_{3} &\longmapsto G_{+} \\
    \end{align*}
    is surjective.
\end{cor}

\begin{proof}
    The element \(G_{-}\) from Proposition \ref{prop:gleason:minus}
    is in the image \(\varphi(R^{\circ})\) by Proposition \ref{prop:gleason:minus}.
    We conclude by Proposition \ref{prop:gleason:surj}.
\end{proof}

We also have the following more general statement:

\begin{cor}
    \label{cor:residue:multivar}
    Consider the map
    \begin{align*}
        \varphi : R \colonequals K\langle T_{1}^{1 / p^{\infty}}, \ldots, T_{n+2}^{1 / p^{\infty}}\rangle  &\longrightarrow 
        K_{r_{1}, \ldots, r_{n}}^{\perfd} \\
        T_{1} &\longmapsto x_{1} \\
        \vdots \\
        T_{n} &\longmapsto x_{n} \\
        T_{n+1} &\longmapsto cx_{1}^{-1}\ldots x_{n}^{-1} \\
        T_{n+2} &\longmapsto G, \\
   \end{align*}
   where \(c \in K\) is a constant such that
    \(s < |cx_{1}^{-1}\ldots x_{n}^{-1}| < 1\),
    and \(G\) is the Gleason element from Proposition \ref{prop:gleason:multivar}
    associated to the lattice basis
    \((1, 0, \ldots, 0), \ldots, (0, \ldots, 0, 1), (-1, \ldots, -1)\).
    Then \(\varphi\) is surjective.
\end{cor}

\begin{proof}
    All the \(W_{i}\) for \(G\) are in the image of \(R^{\circ}\) 
    by taking (algebraic) combinations and \(p\)-th roots.
    We conclude by Proposition \ref{prop:gleason:surj}.
\end{proof}
\begin{proof}[Proof of Theorem \ref{thm:multivar:gleason}]Let $(r_1,\dots, r_l)\in (0, 1]^{l}$ be a free polyradius of length $l\leq n$. By Corollary \ref{cor:residue:multivar}, $K_{r_1,\dots, r_l}^{\perfd}$ is a quotient field of the perfectoid Tate algebra in $l+2$ variables. But the perfectoid Tate algebra in $l+2$ variables is obviously a quotient of the perfectoid Tate algebra in $n+2$ variables.\end{proof} 

\subsection{Noether Normalization}

\begin{thm}
    [Noether normalization for completed perfections, 
    see \cite{Gleason22N}, Proposition 4.5]
    \label{thm:perfection:noether}
    Let $K$ be a perfectoid field of characteristic $p$ and let \(X_{\infty} = \Spa(A_{\infty},A^{+}_{\infty})\) 
    be the completed perfectoid of an affinoid rigid-analytic variety
    \(\Spa(A,A^{+})\).
    Let 
    \(\psi \colon X_{\infty} \to \mathbb{B}_{K}^{n, \perfd}\) 
    be an \'etale map.
    Then we can find a finite etale morphism 
    \(X_{\infty} \xrightarrow{} \mathbb{B}_{K}^{n, \perfd}\).
\end{thm}

\begin{proof}
    Since the \'etale site is stable under perfection,
    any such \'etale map \(\psi\) must 
    be a completed perfection of a map 
    \(X \xrightarrow{} \mathbb{B}^{n}\) 
    at the noetherian level.
    Then, apply \cite[Corollary~D.5]{Zavyalov21-22}
    to get a finite \'etale map at the noetherian
    level, and by the stability of the \'etale
    site again, this induces a finite \'etale
    map at the perfectoid level.
\end{proof}

\begin{cor}
    \label{cor:perftate:finext}
    Let \(T_{n+1,K}^{\perfd}
    = K\langle T_{1}^{1 / p^{\infty}},
    \ldots, T_{n+1}^{1 / p^{\infty}}\rangle \)
    be the perfectoid Tate algebra in \(n+1\) 
    variables.
    Then there exists a quotient field \(L\) of 
    \(T_{n+1,K}^{\perfd}\) such that 
    \(L\) has \(K_{r_{1}, \ldots, r_{n}}^{\perfd}\) 
    as a finite separable extension.
\end{cor}

\begin{proof}
    The map \(T_{n+2,K}^{\perfd} \surj K_{r_{1}, \ldots, r_{n}}\) defined in 
    Corollary \ref{cor:residue:multivar}
    contains the ideal \((T_{1}\cdots T_{n} - c)\) 
    in its kernel.
    Thus, it factors through the
    Zariski closed quotient
    \(A_{\infty} = 
    T_{n+2,K}^{\perfd} /(T_{1}\cdots T_{n} - c)_{\spc}\),
    where the subscript $\spc$ denotes the \textit{spectral radical} (see \cite{Dine22}).
    Since $I$ is an ideal in a perfectoid Tate ring of characteristic $p$,
    \(I_{\spc}\) is obtained by including all 
    \(p\)-th roots of elements in \(I\) 
    and then taking the (topological) closure, see \cite{Dine22}, Remark 2.21.
    The spectral radical is also called the
    \textit{dagger nilradical} in the literature
    (see \cite{camargo-2024-analytic-de-rham}). We see that \(A_{\infty}\) is the completed 
    perfection of the strictly affinoid algebra
    \(A = T_{n+2,K} / (T_{1} \cdots T_{n} - c)\), which is of Krull dimension $n+1$.

    Now, we claim that there exists
    an open immersion 
    \(\Spa(A,A^{\circ}) \xrightarrow{} \mathbb{B}^{n+1}\).
    Indeed, 
    \(A_{\infty}\) is the rational localization of
    \(T_{n+2,K}^{\perfd}\) corresponding
    to the condition \(\{c \leq |T_{1}T_{2}\cdot \ldots \cdot T_{n}|\}\).
    By \cite{Huber1}, Lemma 3.10 (and the fact that perfection preserves the adic spectrum and its
    rational subsets), this implies that $A$ is a rational localization of $T_{n+2, K}$. 
    
    Moreover, we can observe that 
    \(\Spa(A,A^{\circ}) = S_{1} \times_{K} \mathbb{B}^{1}\),
    where 
    \begin{equation*}S_{1} = \{c \leq |T_{1}T_{2} \cdot \ldots \cdot T_{n}|\} 
    \ins \mathbb{B}^{n}.\end{equation*}

    Thus, by Theorem \ref{thm:perfection:noether},
    there exists a finite \'etale map
    \(T_{n+1,K}^{\perfd} \xrightarrow{} A_{\infty}\).
    Let \(L\) be the image of the composite map \(T_{n+1,K}^{\perfd}\to 
    A_{\infty} \xrightarrow{} K_{r_{1}, \ldots, r_{n}}^{\perfd}\).
    Then \(L\) has \(K_{r_{1}, \ldots, r_{n}}^{\perfd}\) as a finite
    separable extension.
\end{proof}

\begin{cor}\label{cor:perftate:finext2}
    The field \(L\) in Corollary \ref{cor:perftate:finext}
    has the following property:
    There exists some polyradius \(r_{1}^{\prime}, \ldots, r_{n}^{\prime}\)
    such that \(L\) is a semi-immediate extension 
    of \(K_{r_{1}^{\prime}, \ldots, r_{n}^{\prime}}^{\perfd}\)
\end{cor}

\begin{proof}
    As \(L\) is the surjective image of a perfectoid Tate algebra,
    it is perfectly topologically finitely generated.
    By Corollary \ref{Temkin's remark 4}, there exists a polyradius $(r_{1}^{\prime},\dots, r_{l}^{\prime})$ such
    that $L$ is a semi-immediate extension of $K_{r_{1}^{\prime},\dots, r_{l}^{\prime}}^{\perfd}$, so
    it remains to prove
    that $l=n$. Since $L$ has the finite extension $K_{r_{1},\dots, r_{n}}^{\perfd}/L$, we have
    \begin{equation*}d_{K}(L)=d_{K}(K_{r_1,\dots, r_n}^{\perfd})=n.\end{equation*}Since
    $L/K_{r_{1}^{\prime},\dots, r_{l}^{\prime}}^{\perfd}$ is semi-immediate, we must then have
    \begin{equation*}n=d_{K}(L)=d_{K}(K_{r_{1}^{\prime},\dots,r_{l}^{\prime}}^{\perfd})=l,\end{equation*}
    as desired.
\end{proof}

\begin{rmk}
    Let \(\mathfrak{m}\) be the kernel of the map
    \(T_{n+1,K}^{\perfd} \surj L\).
    Note that by Theorem \ref{Main theorem, in the body of the paper},
    we conclude that \(\mathfrak{m} \cap T_{n+1,K}^{\coperf}\) 
    is \((0)\).
    Note that we can see this more elementarily 
    by noting that the noether normalization map is injective
    (\'etale and thus dominant).
\end{rmk}



\bibliographystyle{plain}
\bibliography{Bib}

\textsc{Department of Mathematics, University of California San Diego, La Jolla, CA 92093, United States} \newline 

E-mail address: \textsf{ddine@ucsd.edu}
\newline

\textsc{Department of Mathematics, University of California San Diego, La Jolla, CA 92093, United States} \newline

E-mail address: \textsf{jgarzell@ucsd.edu} 

\end{document}